\def\obs#1{{\bf (*** #1 ***)} }

 \def\obs#1{}     % Remova esta linha para rodar a versao 1
\documentclass[twoside,letterpaper,draft,11pt]{amsart}
\usepackage{amsmath,amsthm,latexsym,xspace,amscd,amssymb,mathtools,color,enumerate}

\newcommand\setItemnumber[1]{\setcounter{enumi}{\numexpr#1-1\relax}}

\usepackage[mathscr]{eucal}
\usepackage{amssymb}
\usepackage[all]{xy}
\makeatletter\renewcommand\theenumi{\@roman\c@enumi}\makeatother

\newtheorem{teo1}{Theorem}[section]

\newtheorem{lem1}[teo1]{Lemma}
\newtheorem{cor1}[teo1]{Corollary}
\newtheorem{prop1}[teo1]{Proposition}

\newtheorem{exe}[teo1]{Example}
\newtheorem{remark}[teo1]{Remark}

\newcommand{\m}{{}^{-1}}

\newcommand{\vu}{\vspace{.1cm}}
\newcommand{\vd}{\vspace{.2cm}}

\newcommand{\nod}{\noindent}

\newcommand{\I}{{I}}

\newcommand{\A}{{A}}
\newcommand{\D}{{\mathcal D}}
\newcommand{\tm}[1]{\,\,\text{#1}\,\,}

\newcommand{\af}{\alpha}
\newcommand{\bt}{\beta}
\newcommand{\lb}{\lambda}

\newcommand{\gm}{\gamma}

\newcommand{\ta}{\tau}

\newcommand{\aft}{\af^\ast}

\newcommand{\somai}{\sum_{1\leq i\leq r}}
\newcommand{\somay}{\sum_{y\in \G_0}}
\newcommand{\somag}{\sum_{g\in \G}}
\newcommand{\somaz}{\sum_{z\in\G_0}}
\newcommand{\somal}{\sum_{l\in\G(x)}}

\newcommand{\G}{\mathcal{G}}

\def\ndv{\ {\mid \kern -0.7 em {\scriptstyle \not}} \ \ }

\def\nd{\ {\mid \kern -0.4 em {\scriptstyle \not}} \ \ }

\begin{document}

\thispagestyle{empty}

\title[Lifting partial actions: from groups to groupoids]{Lifting partial actions: from groups to groupoids}

\author[D. Bagio]{Dirceu Bagio}
\address{ Departamento de Matem\'atica, Universidade Federal de Santa Maria, 97105-900\\
	Santa Maria-RS, Brasil}
\email{ bagio@smail.ufsm.br}

\author[A. Paques]{Antonio Paques }
\address{Instituto de Matem\'atica e Estat\'istica, Universidade federal de Porto Alegre, 91509-900\\
	Porto Alegre-RS, Brazil}
\email{paques@mat.ufrgs.br}

\author[H. Pinedo]{H\'ector Pinedo}
\address{Escuela de Matematicas, Universidad Industrial de Santander, Cra. 27 Calle 9 ´ UIS
	Edificio 45, Bucaramanga, Colombia}
\email{ hpinedot@uis.edu.co}


\date{\today}
\begin{abstract} In this paper we are interested in the study of the existence of connections between partial groupoid actions and partial group actions. Precisely, we prove that there exists a datum connecting a partial action of a connected groupoid and a partial action of any of its isotropy groups. 
Furthermore, it will be proved that under a suitable condition the partial skew groupoid ring corresponding to a partial action by a connected groupoid is isomorphic to a specific partial skew group ring.
We also present a Morita theory and a Galois theory related to these partial actions as well as considerations about the strictness of the corresponding Morita contexts.  Semisimplicity, separability and Frobenius properties of the corresponding partial skew groupoid rings are also considered.
\end{abstract}

\maketitle

\setcounter{tocdepth}{1}

\tableofcontents

\section{Introduction}
This paper is motivated by the following three  well known facts:
\begin{enumerate}
\item any groupoid is a disjoint union of its connected components;
\item any groupoid action is uniquely determined by the respective actions of its connected components;
\item any connected groupoid is equivalent to a group as categories.
\end{enumerate}

\vu

Therefore it is natural to expect some kind of connection between groupoid actions and group actions and it is enough to investigate this in the connected case.
Our aim in this paper is to investigate the connections between such both concepts in the most general setting of partial actions.

Throughout this work,  by ring we mean an associative and not necessarily unital ring.
%For instance, it is natural to ask wether both the categories of partial groupoid  actions and the one of the partial group actions are equivalent. As we will see along this paper this is not true in general. Nevertheless, we will present a characterization of the exact conditions under which such an equivalence holds.

\vu

According to \cite{BP} a partial action of a groupoid $\G$ on a ring $A$ is a set $\af=\big(A_g, \af_g\big)_{g\in\G}$ where, for all $g\in\G$, $A_g$ is an ideal of $A_{t(g)}$ ($t(g)=gg\m$), $A_{t(g)}$ is an ideal of $A$ and $\af_g:A_{g\m}\to A_g$ is a ring isomorphism,  for which some appropriate conditions of compatibility hold.

In the case that $\G$ is connected $\af$ gives rise to a specific  datum. To construct such a  datum we start by fixing an object $x$ of $\G$ and a set $\ta(x)$ whose elements are obtained picking up one unique morphism $\ta_y:x\to y$ for each object $y$ of $\G$ (notice that $\ta_y$ always exists for $\G$ is connected).  The  datum arisen from $\af$ is then the triple $\big(A_\af, \af_{\ta(x)}, \af_{(x)}\big)$ where $A_\af$ denotes the set of the ideals $A_y$ indexed by the objects of $\G$, $\af_{\ta(x)}$ denotes the set of the ring isomorphisms $\af_{\ta_y}:A_{\ta\m_y}\to A_{\ta_y}$ for all $\ta_y\in\ta(x)$ and $\af_{(x)}$ denotes the partial action of the isotropy group $\G(x)$ of $\G$ associated to the object $x$ on the ideal $A_x$ of $A$. This correspondence $\af\mapsto \big(A_\af, \af_{\ta(x)}, \af_{(x)}\big)$  gives rise to  a functor from the category of all partial actions of $\G$ on $A$ to the category of all  data of this type.

\vu

Our main purpose in this paper is indeed to show how to construct  partial actions of $\G$ on $A$   coming from  datum of the type above described, the so called \emph{lifted partial groupoid actions}.

\vu

In the next section we recall the respective formal definitions of (connected) groupoid and partial groupoid action, and detailing the fact that every groupoid is a disjoint union of its connected components.

\vu

From the section 3 on we restrict our study to the partial actions by connected groupoids.

\vu

In section 3 we introduce the category $\G_{par}(A)$ of all partial actions of $\G$ on $A$, the category  $\D_{\ta(x)}(\A)$, whose objects are datum of the type
above described, and establish the right conditions for the existence of a good correlation between them via appropriate functors. As consequence we show that any partial groupoid action is an extension of a lifted one (see Proposition \ref{sub} (i)). We also prove that the category $\D_{\ta(x)}(\A)$ is independent of the choice of the object $x$ of $\in\G$ as well as of its transversal $\ta(x)$ (see Theorem \ref{independent}). We exploit  exhaustively such a correlation in the ensuing sections.

\vu

In section 4 we discuss the necessary and sufficient conditions for a lifted partial groupoid  action to be globalizable. In particular we introduce the notion of a \emph{globalizable datum} and describe in details the globalization of the corresponding lifted partial groupoid action (see Theorem \ref{teo:globaliza1}).

\vu

Again according to \cite{BP} any unital partial action $\af$ of a groupoid $\G$ on a ring $A$ gives rise, as expected,  to two new rings, namely, the ring $A^\af$ of the invariants of $A$ under $\af$ and the partial skew groupoid ring $A\star_\af\G$, as well as a canonical Morita context connecting them. Such a context keeps a close relation with the notion of Galois extension. In fact, the condition to ensure that $A$ is a Galois extension of $A^\af$ depends on the strictness of this context. Actually,  by \cite{BP} Theorem 5.3, such a strictness is equivalent to $A$ being a Galois of $A^\af$ and a specific map, called \emph{trace}, from $A$ to $A^\af$ being surjective.  Sections 5, 6 and 7 are dedicated to dealing with all these mentioned concepts.
% in the specific case of lifted partial groupoid actions.

\vu

Section 5 concerns to the partial skew groupoid ring $A\star_\af\G$. We prove, under an appropriate condition, the existence of a ring isomorphism between $A\star_\af\G$ and the partial skew group ring $(A\star_{\af^{\!\ast}}\G_0^2)\star_\theta\G(x)$, where $\G_0$ denotes the set of the objects of $\G$, $\G_0^2=\G_0\times\G_0$ is a groupoid whose structure is standard and well known, $\af^{\!\ast}$ is a global action of $\G_0^2$ on $A$ and $\theta$ is a partial group action of $\G(x)$ on $A\star_{\af^{\!\ast}}\G_0^2$ (see Theorem \ref{la grande finale}). Also, we specialize such condition to the case of lifted partial groupoid actions (see Proposition \ref{58}).

\vu

In sections 6, 7 and 8 we restrict our study only to the lifted partial groupoid actions. In section 6 we recall the notions of invariants and trace map and discuss their properties and correlations (see Theorem \ref{t54}).

\vu

Section 7 is devoted to the Morita theory related to lifted partial groupoid actions.  More specifically, we construct Morita contexts connecting rings of invariants and partial skew groupoid rings and make considerations on their strictness. In particular we show how closed are the corresponding Morita and Galois theories (see Theorem \ref{76}).

\vu

Finally we end this manuscript with the section 8, where we deal with  ring theoretics properties such  the semisimplicity, separability and Frobenius  of the partial skew groupoid ring corresponding to a lifted partial groupoid action (see Theorems \ref{82} and \ref{84}, and Corollary \ref{83}).

\vu

\section{Preliminaries}
As preliminaries we firstly recall from the literature the notions of groupoid, connected groupoid and partial groupoid action and, in the sequel, we show how a groupoid $\G$ can be decomposed in a disjoint union of connected ones. Via such a decomposition it becomes  clear that any partial action of $\G$ is univocally determined by the partial actions of its connected components, which reduces the study of partial groupoid actions to the connected case.

\subsection{Groupoids}
A {\it groupoid} $\G$ is a small category where every morphism is an isomorphism. As a small category $\G$ is composed by a set of morphisms
and  a set  of objects.
Furthermore, to each morphism $h$ of $\G$ correspond naturally the objects $s(h)$ and $t(h)$ of $\G$ called the {\it source} (or domain) and the {\it target} (or range) of $h$ respectively. In particular each object $x$ of $\G$ can be identified with its identity morphism, here also denoted by $x$. Hence, any element of $\G$ can be seen as a morphism of $\G$ and all such considerations can be summarized in the following presentation  $$\G\overset{s, t}\rightrightarrows \G_0$$  where $\G_0$ denotes the set of the objects of $\G$ and $s$ and $t$ the source and the target maps.

\vu

The composition in $\G$ is the map $$m:\G_s\times\!_t\G\to\G \quad\text{denoted by}\quad m(g,h)=gh,$$ where $$\G_s\times\!_t\G=\{(g,h)\in\G\times\G \ |\ s(g)=t(h)\}.$$

\vu

In particular it follows from this characterization of a groupoid that indeed
\begin{align}
\label{eq1} & s(g)=g^{-1}g=t(g^{-1})\,\,\text{ for all}\,\, g\in\G,&\hspace{1.8cm}&\\
\label{eq2} & s(x)=t(x)=x=x^2\,\,\text{ for all}\,\, x\in\G_0,&\hspace{1.8cm}&\\
\label{eq3} & t(g)g=g=gs(g) \,\,\text{ for all}\,\, g\in\G,&\hspace{1.8cm}&\\
\label{eq4} & s(gh)=s(h)\text{ and } t(gh)=t(g),\,\,\text{ for all}\,\, (g,h)\in\G_s\!\times\!\!\,_t\G.&\hspace{1.8cm}&\
\end{align}
Such above relations will be freely used along the text.

\vu

For each $x\in\G_0$ we set $\mathcal{S}_{x}=\{g\in \mathcal{G}: s(g)=x\}$, $\mathcal{T}_x=\{g\in \mathcal{G}: t(g)=x\}$ and $\G(x)=\mathcal{S}_x\cap\mathcal{T}_x$. In particular $\G(x)$ is a group, called the \emph{isotropy (or principal) group associated to $x$.}

\subsection{Connected groupoids}

Let $\G\stackrel{s,t}\rightrightarrows \G_0$ be as in the previous subsection. Any subgroupoid $\mathcal{H}$ of $\G$ is called \emph{connected} if for every $x,y\in\mathcal{H}_0$ the set $\mathcal{H}(x,y)$ of all morphisms from $x$ to $y$ is not empty. Such a notion of connectedness suggests a natural way to obtain a partition of $\G$ in connected components via the following equivalence relation on $\G_0$: for any $x,y\in\G_0$

\begin{equation*}\label{relation}
x\sim y \Leftrightarrow \G(x,y)\neq \emptyset,
\end{equation*}
that is, if and only if there exists $g\in\G$ such that $s(g)=x$ and $t(g)=y$.

Every equivalence class $X\in\G_0/\!\!\sim$ determines a full connected subgroupoid $\G_X \stackrel{s,t}\rightrightarrows X$ of $\G$, which is also called a \emph{connected component} of $\G$. It is immediate to see that $$\G=\dot\cup_{X\in \G_0/\!\sim}\G_X$$ (where the symbol $\dot\cup$ denotes a disjoint union).

\subsection{Partial groupoid action}
According to \cite{BP}, a \emph{partial action} of a groupoid $\G$ on  a ring $\A$ is a pair
$\af=(\A_g,\af_g)_{g\in \G}$ where, for each $g\in \G,$  $\A_{t(g)}$ is an ideal of $\A$, $\A_g$ is an ideal of $\A_{t(g)}$, $\af_g:\A_{g\m}\to \A_g$ is an isomorphism of rings,
and the following conditions hold:
\begin{align}
\label{par1} &  \af_x\tm{is the identity map} \tm{of} \A_x,& &\hspace{5.5cm}&\\
\label{par2} & \af_g\af_h\le \af_{gh},& &\hspace{5.5cm}&
\end{align}
for all $x\in \G_0$ and $(g,h)\in\G_s\!\times\!\!\,_t\G.$
We say that $\af$ is {\it global} if $\af_g\af_h=\af_{gh}$, for all $(g,h)\in \G_s\!\times\!\!\,_t\G.$

\vu

The condition \eqref{par2} above means that $\af_{gh}$ is an extension of $\af_g\af_h$, that is, the domain $\A_{(gh)^{-1}}$ of $\af_{gh}$ contains the domain $\af_h^{-1}(\A_{g^{-1}}\cap\A_h)$ of $\af_g\af_h$, and
both the maps coincide on this last set.

\vu

Notice also that $\af$ induces by restriction a partial action $\af_{(x)}=(\A_g,\af_g)_{g\in\G(x)}$ of the group $\G(x)$ on the ring $\A_x$, for all $x\in\G_0$.

\vu

Lemma 1.1 from \cite{BP}  gives some complementary properties of partial groupoid actions, enumerated below, that will also be useful in the sequel.
\begin{align}
\label{bapa1} &\af\tm{is global if and only if} \A_g=\A_{t(g)},\tm{for all} g\in \G,& &\hspace{2.5cm}&\\
\label{bapa2} &\af_{g\m}=\af\m_g,\tm{for all} g\in \G,& &\hspace{2.5cm}&\\
\label{bapa3}&\af_g(\A_{g\m}\cap \A_h)=\A_{g}\cap \A_{gh},\tm{for all} (g,h)\in \G_s\!\times\! _t\G.& &\hspace{2.5cm}&\
\end{align}

\begin{remark} {\rm It is clear that given any two disjoint  groupoids $\G'\stackrel{s,t}\rightrightarrows \G'_0$ and $ \G''\stackrel{s,t}\rightrightarrows \G''_0$ one always can form the groupoid $\G=\G'\cup \G'' \stackrel{s,t}\rightrightarrows\G_0=\G'_0\cup \G''_0,$ whose composition is induced in an obvious way by the compositions in $\G'$ and $\G''$ respectively. Furthermore, it is also straightforward to check that partial actions of $\G$ on  a ring $\A$ induce by restriction partial actions of $\G'$ and $\G''$ on $\A$ and, conversely, they are univocally determined by the partial actions of $\G'$ and $\G''$ on $\A$. Therefore, it follows from the subsection 2.2 that the study of partial groupoid actions reduces to the connected case.}
\end{remark}

\section{Lifting partial group actions}

From now on by $\G\stackrel{s,t}\rightrightarrows \G_0$ we will always denote a connected groupoid and by $\A$  a ring on which $\G$ acts partially. Recalling notations, $$\G(x,y)=\{g\in\G\ |\ s(g)=x\,\ \text{and}\,\ t(g)=y\}=\mathcal{S}_x\cap\mathcal{T}_y$$ and, in particular, $$\G(x,x)=\mathcal{S}_x\cap\mathcal{T}_x=\G(x),$$ for all $x,y\in\G_0$. Furthermore, recall that the identity map of each object $y\in\G_0$ is also denoted by $y$.

\vu

Hereafter, let $x\in \G_0$ be a fixed object of $\G$. Now consider on $\mathcal{S}_x$ the following equivalence relation:
$g\equiv_x l \Leftrightarrow l\m g\in\G(x)\Leftrightarrow t(g)=t(l)$. A transversal for $\equiv_x$  such that $\tau_x=x$ will be called a \emph{transversal for $x$} and denoted by $\ta(x)$, that is,
$\tau(x)=\{\tau_{y}:y\in \G_0\}$ where $\ta_y$ is a chosen morphism in $\G(x,y)$, for each $y\in\G_0$ with  $\ta_x=x$.

\vu

Observe that any  tranversal $\ta(x)$ determines a natural map $\pi$ from $\G$ to $\G(x)$ that associates to each element  $g\in\G$ a unique element $g_x\in\G(x)$ via the following formula
\begin{align}\label{epi}
\pi(g)=g_x=\ta\m_{t(g)}g\ta_{s(g)},
\end{align}
illustrated by the following diagram
\[\xymatrix{  s(g) \ar[rr]^{g} & & \ar[ld]^{\tau\m_{t(g)}} t(g) \\
	          &x\ar[lu]^{\tau_{s(g)}} & }\]

\vu

\begin{lem1}\label{lema31}
$\pi$ is a	groupoid epimorphism such that
\begin{align}\label{epicor}
\pi(\ta(x))=\{x\}\,\text{ and }\, \pi(h)=h,\,\text{ for all }\,h\in \G(x).
\end{align}

\end{lem1}
\begin{proof} Observe that $\pi(\G_0)=\{x\}$. Indeed, $$y_x=\ta\m_{t(y)}y\ta_{s(y)}\overset{\eqref{eq2}}{=}\ta\m_yy\ta_y=\ta\m_y\ta_y=x.$$
Moreover, for all $g,h\in \G$ such that $s(g)=t(h)$, we have
\begin{align*}(gh)_x&=\tau\m_{t(gh)}gh\tau_{s(gh)}\\
&\overset{\mathclap{\eqref{eq4}}}{=}\tau\m_{t(g)}gh\tau_{s(h)}\\
%&\overset{\mathclap{\eqref{eq3}}}{=}\tau\m_{t(g)}gt(h)h\tau_{s(h)}\\
&=\ta\m_{t(g)}g(\ta_{t(h)}\ta\m_{t(h)})h\ta_{s(h)}\\
&=(\ta\m_{t(g)}g\ta_{s(g)})(\ta\m_{t(h)}h\ta_{s(h)})\\
&=g_xh_x.\end{align*}

Also, for all $y\in\G_0$,
\begin{align*}
(\ta_y)_x=\ta\m_{t(\ta_y)}\ta_y\ta_{s(\ta_y)}=\ta\m_y\ta_y\ta_x\overset{\eqref{eq1}}{=}s(\ta_y)\ta_x=x\ta_x=x^{2}\overset{\eqref{eq2}}{=}x.
\end{align*}

Finally, if $h\in\G(x)$ then
$$h_x=\ta\m_{t(h)}h\ta_{s(h)}=\ta\m_xh\ta_x=xhx=h.$$
\end{proof}

\subsection{The category $\G_{par}(\A)$}

We will denote by $\G_{par}(\A)$ the category  whose objects are  partial actions of $\G$ on  a fixed ring $\A$ and whose morphisms are defined as follows. Given  $\af=(\A_g,\af_g)_{g\in \G}$  and $\af'=(\A'_g,\af'_g)_{g\in \G}$ in  $\G_{par}(\A)$, a morphism $\psi\colon \af\to \af'$ is a set of ring homomorphisms
$\psi=\{\psi_{y}:A_{y}\to A'_{y}\}_{y\in \G_0}$ such that

\begin{enumerate}[$\quad$(a)]
\item $\psi_{t(g)}( \A_g)\subseteq \A'_g,$

\vu

\item $\af'_g\psi_{s(g)}=\psi_{t(g)}\af_g\,\,\text{ in}\,\, \A_{g\m},$
\end{enumerate}
for all $g\in \G.$

\vd

Observe now that any pair $(\af,\ta(x))$, where $\af=(\A_g,\af_g)_{g\in \G}\in\G_{par}(\A)$ and $\ta(x)$ is a transversal for $x$,  determines the following  datum:
\begin{enumerate}[$\quad\circ $]
\item  a set of ideals of $\A$: $\A_\af=\{\A_y\}_{y\in\G_0}$,

\vu

\item  a set of ring isomorphisms $\af_{\ta(x)}=\{\af_{\ta_{y}}\colon \A_{\ta\m_{y}}\to \A_{\ta_{y}}\}_{y\in \G_0}$ and

\vu

\item  a partial group action $\af_{(x)}=(\A_g, \af_g)_{g\in \G(x)}$ of $\G(x)$ on $\A_x$.
\end{enumerate}
Moreover, since $t(\ta_{y})=y$ (resp., $t(\ta^{-1}_{y})=x$), then $\A_{\ta_{y}}$ (resp., $\A_{\ta^{-1}_{y}}$) is an ideal of $\A_y$ (resp., $\A_x$), for all $y\in\G_0$. In particular, $\A_{\ta^{-1}_{x}}=\A_{\ta_{x}}=\A_x$ and $\af_{\ta_x}=\af_x$ is the identity map of $A_x$ by \eqref{par1}.

\vu

Finally, observe that
\begin{align*}
\af_{\ta_{t(g)}}(A_{\ta\m_{t(g)}}\cap \af_{g_x}(A_{\ta\m_{s(g)}}\cap A_{g\m_x}))&\overset{\mathclap{\eqref{bapa3}}}{=}\af_{\ta_{t(g)}}(A_{\ta\m_{t(g)}}\cap A_{\ta\m_{t(g)}g}\cap A_{g_x})\\
&\overset{\mathclap{\eqref{bapa3}}}{=}A_{\ta_{t(g)}}\cap A_{g}\cap A_{g\ta_{s(g)}},
\end{align*}
which is an ideal of $A_{t(g)}$, for all $g\in \G$.

\vd

This suggests the existence of a new category, named  here a \emph{data category}, which we will describe below.

\subsection{The category $\mathcal{D}_{\ta(x)}(\A)$}

\vu

To describe the  data category $\D_{\ta(x)}(\A)$, with $\ta(x)$ a fixed transversal, our inspiration is the  datum we point out in the previous subsection.
The objects of $\D_{\ta(x)}(\A)$ are given by  datum of the following type:
\begin{enumerate}[$\quad$(a)]
\item  a set of ideals of $\A$: $I=\{I_y\}_{y\in\G_0}$,

\vu

\item  a set of ring isomorphisms: $\gamma_{\ta(x)}=\{\gamma_{\ta_{y}}:\I_{\ta_{y}^{-1}}\to\I_{\ta_{y}}\}_{y\in\G_0}$ and

\vu

\item  a partial group action $\gm_{(x)}=(\I_g,\gm_g)_{g\in\G(x)}$ of $\G(x)$ on $\I_x$,
\end{enumerate}
under the additional assumptions
\begin{enumerate}[$\quad$(a)]
	  \setItemnumber{4}
	
\item $\gamma_x=\gamma_{\ta_x}$ is the identity map of $I_x=I_{\ta_x}$,
	
	\vu
	
\item $\I_{\ta_{y}}$ (resp., $\I_{\ta_{y}^{-1}}$) is an ideal of $\I_{t(\ta_{y})}=\I_y$ (resp., $\I_{t(\ta_{y}^{-1})}=\I_x$), for all $y\in \G_0$,
	
	\vu
	
\item   $\gamma_{\ta_{t(g)}}(I_{\ta\m_{t(g)}}\cap \gamma_{g_x}(I_{\ta\m_{s(g)}}\cap I_{g\m_x}))$ is an ideal of $I_{t(g)}$, for all $g\in \G$.
\end{enumerate}

\vu

For each such a  datum we denote the corresponding object of $\D_{\ta(x)}(\A)$ by the triple $(\I, \gamma_{\ta(x)},\gm_{(x)})$.
The morphisms of $\D_{\ta(x)}(\A)$ are given by sets of ring homomorphisms denoted and described as follows:
$$f=\{f_y\}_{y\in \G_0}\colon (\I, \gamma_{\ta(x)}, \gm_{(x)})\to (\I', \gamma'_{\ta(x)}, \gm'_{(x)})$$ where
\begin{enumerate}[$\quad\circ $]
\item  for each $y\in \G_0,$ $f_y\colon \I_y\to \I'_y$  is a ring homomorphism such that $$f_y(\I_{\ta_{y}} )\subseteq \I'_{\ta_{y}}\,\,(\text{if}\,\  y\neq x),\,\,  f_x(\I_{\ta\m_{y}} )\subseteq \I'_{\ta\m_{y}}\,\ \text{ and}\,\, \,\gamma'_{\ta_{y}}f_x=f_y\gamma_{\ta_{y}}.$$

\vu

\item $f_x\colon \gm_{(x)} \to \gm'_{(x)}$ is a morphism of partial group actions, i.~e, $f_x\colon \I_x\to \I'_x$ satisfies

\[f_x(\I_h)\subseteq \I'_h,\qquad \gm'_hf_x=f_x\gm_h \ \text{ in } \ \I_{h\m},\ \text{ for all } \ h\in \G(x).\]
\end{enumerate}

\vu

Now, we prove that the category $\D_{\ta(x)}(A)$ does not depend neither on the choice of the object $x\in\G_0$, nor on the choice of the transversal $\ta(x)$ of $x$.

\begin{prop1}\label{independent}
	$\D_{\ta(x)}(A)$ and $\D_{\lb(z)}(A)$ are isomorphic as categories, for all $x, z\in\G_0$ and their respective transversals $\ta(x)$ and $\lb(z)$.
\end{prop1}

\begin{proof}
	We start by defining the functor $F^{\lb(z)}_{\ta(x)}:\D_{\ta(x)}(A)\to \D_{\lb(z)}(A)$, its reverse is defined similarly.
	
	Firstly the correspondence between the objects:  to each $\gamma=(I,\gamma_{\ta(x)},\gamma_{(x)})\in\D_{\ta(x)}(A)$ we associate $\gamma'=(I',\gamma'_{\ta(x)},\gamma'_{(z)})\in\D_{\lb(z)}(A)$ chosen in the following way:
	\begin{enumerate}[$\quad\circ $]
		\item  if $I=\{I_y\}_{y\in\G_0}$ we take $I'=\{I'_y\}_{y\in\G_0}$ such that $I'_z=I_x$, $I'_x=I_z$ and $I'_y=I_y$ for all $y\notin\{ x, z\}$,
		\item  if $\gamma_{\ta(x)}=\{ \gamma_{\ta_{y}}:I_{\ta\m_y}\to I_{\ta_y}\}_{y\in\G_0}$ we take $\gamma'_{\lb(z)}=\{ \gamma'_{\lb_{z}}:I'_{\lb\m_z}\to I'_{\lb_z}\}_{z\in\G_0}$
		such that $\gamma'_z=\gamma'_{\lb_z}=\gamma_{\ta_{(x)}}=\gamma_x$, $\gamma'_{\lb_x}=\gamma_{\ta_z}$ and $\gamma'_{\lb_y}=\gamma_{\ta_y}$
		for all $y\notin\{ z, x\}$,
		\item  if $\gamma_{(x)}=(I_h,\gamma_h)_{h\in\G(x)}$ we take  $\gamma'_{(x)}=(I'_l,\gamma'_l)_{l\in\G(z)}$ such that $I'_l=I_{\phi(l)}$, $\gamma'_l=\gamma_{\phi(l)}$ for
		all $l\in\G(z)$, where $\phi:\G(z)\to\G(x)$ is the group isomorphism defined by $l\mapsto \ta\m_zl\ta_z$.
	\end{enumerate}
	
	Now the correspondence between the morphisms: for a morphism $(f_y, f^{(x)}):\gamma\to\delta$ in $\D_{\ta(x)}(A)$ we associate the morphism $(f'_y, f'^{(z)}):\gamma'\to\delta'$ in $\D_{\lb(z)}(A)$,
	given by $f'^{(z)}=f^{(x)}$, $f'_x=f_z$ and $f'_y=f_y$ for all $y\notin\{ x, z\}$.
	
	This way we have got the functors $F_{\ta(x)}^{\lb(z)}$ and $F_{\lb(z)}^{\ta(x)}$. It is straightforward to check that the compositions
	$F_{\lb(z)}^{\ta(x)}\circ F_{\ta(x)}^{\lb(z)}$ and $F_{\ta(x)}^{\lb(z)}\circ F_{\lb(z)}^{\ta(x)}$ are respectively the identity functors
	of $\D_{\ta(x)}(A)$ and $\D_{\lb(z)}(A)$.
\end{proof}

\subsection{The functors $F_{\ta(x)}$ and $G_{\ta(x)}$}

Our purpose in this subsection is to find  functors relating both the categories $\G_{par}(\A)$ and $\D_{\ta(x)}(\A)$.
It is easy to see from the subsection $3.1$ that the association
$F_{\ta(x)}\colon \G_{par}(\A)\to \D_{\ta(x)}(\A)$ given by
\begin{align*}
&\qquad\qquad\qquad\qquad \af\longmapsto (\A_\af, \af_{\ta(x)}, \af_{(x)}),&\\[.2em]
&\af\overset{\psi}{\longrightarrow}\af'\longmapsto F_{\ta(x)}(\af)\overset{F_{\ta(x)}(\psi)}{\longrightarrow} F_{\ta(x)}(\af'), \ \text{ where } F_{\ta(x)}(\psi)=\{\psi_y\}_{y\in \G_0},&
\end{align*}
is indeed a functor.

\vu

For the reverse functor, given a triple $(\I, \gamma_{\ta(x)},\gm_{(x)})$ we set $\bt=(B_g,\bt_g)_{g\in \G}$, where each $\bt_g$ is the following ring isomorphism
\begin{equation}\label{afcon}\bt_g=\begin{cases}
\text{the identity map of } I_{y},& \text{ if } g=y\in \G_0,\\
\gamma_{\ta_{t(g)}}\circ\gamma_{g_x}\circ\gamma\m_{\ta_{s(g)}},& \text{ if } g\notin \G_0.
\end{cases}
\end{equation}
and each $B_g$ is taken as the range of $\bt_g,$ for all $g\in \G.$

\vd

In order to completely understand \eqref{afcon} it is convenient to recall how are defined the domain and the range of a composition of two partial bijections. Notice that the maps in \eqref{afcon} are partial bijections of $A$.
In general, if $u$ and $v$ are two given partial bijections of a non empty set, the domain and the range of the composition $u\circ v$ are respectively defined as follows:
$$\text{dom}(u\circ v)=v^{-1}(\text{dom}(u)\cap \text{ran}(v))\quad\text{and}\quad \text{ran}(u\circ v)= u(\text{dom}(u)\cap\text{ ran}(v)).$$

In the specific case of \eqref{afcon} we have
\begin{align}
\label{dom} B_{g}=\begin{cases}
I_{y},& \text{ if } g=y\in \G_0,\\
\gamma_{\ta_{t(g)}}(I_{\ta\m_{t(g)}}\cap \gamma_{g_x}(I_{\ta\m_{s(g)}}\cap I_{g\m_{x}})),& \text{ if } g\notin \G_0.
\end{cases}
\end{align}

\begin{teo1}\label{outro} The pair $\bt=(B_{g},\bt_{g})_{g\in \G}$ constructed above  is an object in the category $\G_{ par}(\A)$.
\end{teo1}
\begin{proof}
Take  $g, h\notin \G_0$ such that $s(g)=t(h)$. Then, by restriction to the domain of $\bt_g\bt_h$ we have the following inequalities
\begin{align*}
\bt_{g}\bt_{h}&=\gamma_{\ta_{t(g)}}\gamma_{g_x}\gamma\m_{\ta_{s(g)}}\gamma_{\ta_{t(h)}}\gamma_{h_x}\gamma\m_{\ta_{s(h)}}\\
%&=\gamma_{\ta_{t(g)}}\gamma_{g_x}\gamma\m_{\ta_{s(g)}}\gamma_{\ta_{s(g)}}\gamma_{h_x}\gamma\m_{\ta_{s(h)}}\\
&\overset{}{\leq}\gamma_{\ta_{t(g)}}\gamma_{g_x}\gamma_{h_x}\gamma\m_{\ta_{s(h)}}\\
&\overset{\mathclap{(\ast)}}{\le} \gamma_{\ta_{t(g)}}\gamma_{g_xh_x}\gamma\m_{\ta_{s(h)}}\\
&\overset{\mathclap{\eqref{epi}}}{=}\gamma_{\ta_{t(g)}}\gamma_{(gh)_x}\gamma\m_{\ta_{s(h)}}\\
&\overset{\mathclap{\eqref{eq4}}}{=}\gamma_{\ta_{t(gh)}}\gamma_{(gh)_x}\gamma\m_{\ta_{s(gh)}}\\
&=\bt_{gh}.\end{align*}
Notice that $(\ast)$  holds for %$\gamma_{h_x}=\gamma\m_{\ta_{s(g)}}\gamma_{\ta_{s(g)}}\gamma_{h_x}$, and $(\ast)$ 
follows for $\gamma$ is a partial group action of $\G(x)$ over $I_x$. Since $\text{dom}(\bt_g\bt_h)\subseteq\text{dom}(\bt_{gh})$, we obtain $\bt_g\bt_h\leq \bt_{gh}$.

\vu

If $h=y\in \G_0$ (resp. $g=y\in \G_0$ ) and $s(g)=y$ (resp. $t(h)=y$) then $\beta_{g}\bt_h=\beta_g=\beta_{gh}$ (resp. $\beta_{g}\bt_h=\beta_h=\beta_{gh}$). When $g=h=y\in \G_0$, we have $\bt_{g}\bt_{h}=\beta_y=\bt_{gh}$.

\vu

By assumption, $B_y=I_y$ is an ideal of $A$ and $B_g$ is an ideal of $B_{t(g)}=I_{t(g)}$, for all $y\in \G_0$ and $g\notin G_0$.
\end{proof}

\begin{remark}\label{r32}
	{\rm It is clear from \eqref{afcon} and \eqref{dom} that $\bt_g=\gamma_g$ and $B_g=I_g$ for all $g\in\G(x)$, that is, $\bt_{(x)}=\gamma_{(x)}$. Thence, any partial action $\bt=(B_g,\bt_g)_{g\in\G}$ of $\G$ on $A$, constructed as above from a given  datum $\gamma=(\I, \gamma_{\ta(x)},\gm_{(x)})$, will be referred as a \emph{lifted partial groupoid action}}.
\end{remark}

By Theorem \ref{outro}, the association $G_{\ta(x)}:\D_{\ta(x)}(A)\to\G_{par}(A)$ given by
\begin{align*}
&G_{\ta(x)}((\I, \gamma_{\ta(x)},\gm_{(x)}))=\bt=(B_g,\bt_g)_{g\in\G},\qquad \text{(see }  \eqref{afcon} \ \text{ and } \ \eqref{dom}\text{)}&\\[0.2em]
&\qquad\qquad\qquad G_{\ta(x)}(({f_y}, f^{(x)})_{y\in \G_0})=\{f_y\}_{y\in \G_0},&
\end{align*}
is a functor.

\vd

The next result establishes the relation between the functors $F_{\ta(x)}$  and $G_{\ta(x)}$.

\vu

\begin{prop1}\label{sub} Let $F_{\ta(x)}$  and $G_{\ta(x)}$ be the functors constructed above.
\begin{enumerate}[\rm (i)]
	\vu
	
\item If $\af=(A_g,\af_g)_{g\in\G}\in \G_{par}(A)$ and $\beta=G_{\ta{(x)}}\circ F_{\ta{(x)}}(\af)$ then $\beta\leq \af$, that is, $\bt_g\leq \af_g$, for all $g\in \G$. Moreover,   $\beta=\af$ if and only if $A_g\subset A_{\ta_{t(g)}}\cap A_{g\ta_{s(g)}}$, for all $g\notin\G_0$.

\vu

\item  $F_{\ta{(x)}}$ is a left inverse functor of $G_{\ta{(x)}}$, that is, $F_{\ta{(x)}}\circ G_{\ta{(x)}}(\gamma)=\gamma$, for all $\gamma\in\D_{\ta(x)}(A)$.
\end{enumerate}

\vu

\end{prop1}

\begin{proof} Observe that for $g\notin \G_0$ we have
\begin{align*}
\,\,\,\bt_g &\overset{\mathclap{\eqref{afcon}}}{=}\af_{\ta_{t(g)}}\circ\af_{g_x}\circ\af_{\ta\m_{s(g)}}& &\hspace{2cm}& \\
      &\overset{\mathclap{\eqref{epi}}}{=}\af_{\ta_{t(g)}}\circ\af_{\ta\m_{t(g)}g\ta_{s(g)}}\circ\af_{\ta\m_{s(g)}}& &\hspace{2cm}&\\
      &\overset{\mathclap{\eqref{par2}}}{\leq}\af_{\ta_{t(g)}\ta\m_{t(g)}g\ta_{s(g)}\ta\m_{s(g)}}=\af_g.& &\hspace{2cm}&
\end{align*}
If $g=y\in \G_0$, then it is immediate that $\beta_g=\beta_y=\af_y=\af_g$ and the first part of (i) follows. Since
\begin{align*}
B_g &\overset{\mathclap{\eqref{dom}}}{=}\af_{\ta_{t(g)}}(A_{\ta\m_{t(g)}}\cap\af_{g_x}(A_{\ta\m_{s(g)}}\cap A_{g\m_x}))\\
&\overset{\mathclap{\eqref{bapa3}}}{=}\af_{\ta_{t(g)}}(A_{\ta\m_{t(g)}}\cap A_{g_x\ta\m_{s(g)}}\cap A_{g_x})\\
&\overset{\mathclap{\eqref{epi}}}{=}\af_{\ta_{t(g)}}(A_{\ta\m_{t(g)}}\cap A_{\ta\m_{t(g)}g}\cap A_{g_x})\\
&\overset{\mathclap{\eqref{bapa3}}}{=} A_{\ta_{t(g)}}\cap A_g\cap A_{\ta_{t(g)}g_x}\\
&\overset{\mathclap{\eqref{epi}}}{=}A_g\cap(A_{\ta_{t(g)}}\cap A_{g\ta_{s(g)}}),
\end{align*}
for all $g\notin \G_0$, we obtain the second part of (i). The item (ii) is easily checked using the definitions of $G_{\ta{(x)}}$ and $F_{\ta{(x)}}$.
\end{proof}

\vu

A partial action $\af=(A_g,\af_g)_{g\in \G}$ of $\G$ on $A$ will be called {\it $\tau(x)$-global} if \vu
\begin{align}
	\label{cond1} A_{\tau\m_y}=A_x \ \ \text{and} \ \ A_{\tau_y}=A_y, \ \ \text{ for all } \ y\in\G_0.
\end{align}

\vu

Clearly, any global action of $\G$ on $A$ is $\tau(x)$-global. We present below an example of a partial (not global) groupoid action which is $\tau(x)$-global.

\begin{exe}\label{57}
	{\rm Let $\G=\{g,h,l,m,l^{-1},m^{-1}\}\rightrightarrows \{x,y\}=\G_0$ be the  groupoid given in Example \ref{hexd} and
		$A=\mathbb{C}e_1\oplus\mathbb{C}e_2\oplus\mathbb{C}e_3\oplus\mathbb{C}e_4$, where $\mathbb{C}$ denotes the complex number field.
		Consider the partial action $\big(A_z,\af_z\big)_{z\in\G}$ of $\G$ on $A$ described bellow:
		\begin{align*}
			&A_x=\mathbb{C}e_1\oplus\mathbb{C}e_2=A_{l\m},& & A_y=\mathbb{C}e_3\oplus\mathbb{C}e_4=A_l, &\\[.2em]
			& A_g=\mathbb{C}e_1=A_{g\m}=A_{m\m},& & A_m=A_h=\mathbb{C}e_3=A_{h\m}, &\
		\end{align*}
		and
		\begin{align*}
			&\af_x=id_{A_x},\ \ \ \af_y= id_{A_y},\ \  \ \af_g:ae_1\mapsto \overline{a}e_1, \ \ \ \af_h:ae_3\mapsto \overline{a}e_3,\ \  \ \af_m:ae_1\mapsto \overline{a}e_3, \\[.3em]
			&\ \ \af_{m\m}: ae_3\mapsto\overline{a}e_1,\ \ \ \af_l:ae_1+be_2\mapsto ae_3+be_4, \ \ \ \af_{l\m}:ae_3+be_4\mapsto ae_1+be_2,
		\end{align*}
		where $\overline{a}$ denotes the conjugate of $a$, for all $a\in\mathbb{C}$. It is clear that $\alpha$ satisfies the required for the transversal $\ta(x)=\{\ta_x=x, \ta_y=l\}$ of $x$.}
\end{exe}

\begin{remark}
	{\rm The notion of partial action $\ta(x)$-global indeed depends on the choice of the transversal $\tau(x)$. Notice that the partial action $\alpha$ of $\G$ on $A$ given in the previous example is $\tau(x)$-global but not $\lambda(x)$-global, where $\lambda(x)=\{\lambda_x=x, \lambda_y=m\}$.    }
\end{remark}

\vu

The full subcategory of $\D_{\ta(x)}(A)$ with objects  $(\I, \gamma_{\ta(x)},\gm_{(x)})$ satisfying $I_{\ta_y}=I_{y}$ and $I_{\ta\m_y}=I_{x}$, for all $y\in \G_0$, will be denoted by $\D^{gl}_{\ta(x)}(A)$.
Also, the full subcategory of $\G_{par}(A)$  whose objects are the $\ta(x)$-global actions of $\G$ on $A$ will be denoted by 	$\G_{\ta(x)}^{gl}(A)$. 

\vu
\begin{cor1} The categories $\G_{\ta(x)}^{gl}(A)$ and $\D^{gl}_{\ta(x)}(A)$ are isomorphic.
\end{cor1}
\begin{proof}
	Let  $\af=(A_g,\af_g)_{g\in\G}\in \G_{\ta(x)}^{gl}(A)$. Then
	\begin{align*}
	A_g&=\af_g(A_{g\m}\cap A_{{s(g)}})  \\
	&\overset{\mathclap{\eqref{cond1}}}{=}\af_g(A_{g\m}\cap A_{\ta_{s(g)}})\\
	&\overset{\mathclap{\eqref{bapa3}}}{=}A_{g}\cap A_{g\ta_{s(g)}}\\
	&\subset A_{{t(g)}}\cap A_{g\tau_{s(g)}}\\
	&=A_{\tau_{t(g)}}\cap A_{g\tau_{s(g)}},
	\end{align*}
	for all $g\in \G$. It follows from Proposition \ref{sub} (i) that $G_{\ta{(x)}}\circ F_{\ta{(x)}}(\af)=\af$ and whence $F_{\ta{(x)}}$ is a right inverse of $G_{\ta{(x)}}$. Using Proposition \ref{sub} (ii), we conclude the result.
\end{proof}

\vu

We end this subsection with a pair of examples to illustrate the construction of a lifted partial groupoid action.
\begin{exe}\label{hexd} {\rm Let $\G=\{g,h,l,m,l^{-1},m^{-1}\}\rightrightarrows \{x,y\}=\G_0$ the groupoid with the following composition rules
		\[ g^2=x,\quad h^{2}=y,\quad lg=m=hl,\quad g\in\G(x),\quad h\in\G(y)\,\,\text{ and }\,\, l,m\in\G(x,y). \]
		The diagram bellow illustrates the structure of $\G$:
		\[\xymatrix{& x\ar[r]^{l}  &y \ar[d]^{h}\\
			& x\ar[u]^{g} \ar[r]^{m} & y} \]
Fix a transversal $\ta(x)=\{\ta_x=x, \ta_y=l\}$ of $x$,  a ring $A$ and $J,L,I$ ideals of $A$ such that $J,L\subset I$. Fix also ring automorphisms $\gamma$ and $\sigma$ of $I$ such that $\sigma(L)=L$, $\sigma|_L=\sigma\m|_L$ and $\gamma(J)=J$. In all what follows $id_I$ denotes the identity map of $I$.	\vu

Now take the following  datum$\big(I,\gamma_{\ta(x)},\gamma_{(x)}\big)$, where

\vu

\begin{enumerate}[$\quad\circ $]
	\item $\{I_x=I_y=I\}$ is the set of the ideals of $A$ indexed by the elements of $\G_0$,
	
	\vu
	
	\item  $\gamma_{\ta(x)}=\{\gamma_{\ta_x}=\gamma_x=id_{I_x}=id_I, \gamma_{\ta_y}=\gamma_l=\gamma\}$ is the set of the ring isomorphisms,
	
	\vu
	
	\item  $\gamma_{(x)}=\big(\{I_x=I, I_g=L=I_{g\m}\}, \{\gamma_x=id_{I_x}, \gamma_g=\sigma=\gamma_{g\m}\}\big)$ is the partial action of $\G(x)=\{x, g\}$ on $I_x$.
\end{enumerate}	
Applying the functor $G_{\ta{(x)}}$ in the datum$\big(I,\gamma_{\ta(x)},\gamma_{(x)}\big)$ we obtain the following partial action $\bt=(B_u,\bt_u)_{u\in\G}$ of $\G$ on $A$. Firstly, note that \eqref{epi} implies that
$$g_x=h_x=m_x=m\m_x=g,\qquad l_x=l\m_x=x_x=y_x=x.$$
From \eqref{afcon} and \eqref{dom} it follows that
\begin{align*}
&\bt_x=\bt_y=id_I,\qquad \bt_g=\sigma|_{L}=\beta_{g\m},\qquad \bt_h=\gamma\sigma\gamma\m|_{\gamma(J\cap\sigma(J\cap L))}=\beta_{h\m},\\[.5em]
& \bt_m=\gamma\sigma:\sigma(J\cap L)\to \gamma(J\cap L), \quad\bt_{m\m}=\sigma\gamma\m:\gamma(J\cap L)\to \sigma(J\cap L).
\end{align*}}		
\end{exe}

\begin{exe}  {\rm Let $A$ be a unital ring,  $m$ be a positive integer and $\{e_i\}_{i=1}^m$ a set of orthogonal idempotents in the center of $A$  whose sum is $1_A$.  Fix $i_0\in \{1,\cdots, m\}$ and  consider a group  $G$ acting on $Ae_{i_0}$  by ring isomorphisms $\sigma_g$, $g\in G$.  \vu
		
Denote by $\Gamma_G^m$ the groupoid having as objects the set $\{1,\cdots, m\}$ and morphisms $(g,i,j),$ where $g\in G$ and $i,j\in\{1,\cdots, m\}. $ The source and target maps on $\Gamma_G^m$ are  $s(g,i,j)=j$ and $t(g,i,j)=i,$ respectively. The composition   is given by the rule  $(g,i,j)(h,j,k)=(gh,i,k)$. Clearly,  $\Gamma_G^m$ is a connected groupoid. Fix the transversal $\tau(i_0)=\{\tau_i=(e,i,i_0)\}_{1\leq i\leq m}$, where $e$ denotes the identity element of $G$. Consider the following datum $\big(I,\gamma_{\ta(i_0)},\gamma_{(i_0)}\big)$:  }

\vu

\begin{enumerate}[$\quad\circ $]
	\item {\rm $\{I_i=I_{\tau_i}=Ae_{i}\}_{1\leq i\leq m}$ and $\{I_{\tau^{-1}_i}=Ae_{i_0}\}_{1\leq i\leq m} $ is the set of the ideals of $A$ indexed by the objects of $\Gamma_G^m$,}
	
	\vu
	
	\item {\rm $\gamma_{\ta({i_0})}=\{\gamma_{\ta_i}=A{e_{i_0}}\ni ae_i\mapsto ae_{i_0}\in A{e_i}\}$   is the set of the ring isomorphisms,}
	
	\vu
	
	\item  {\rm $\gamma_{(i_0)}=\{\gamma_{{(g,i_0,i_0)}}=\sigma_g\}$ is the action of the group  $\Gamma_G^m(i_0)=\{(g,i_0,i_0)\mid g\in G\}$ on $A{e_{i_0}}$.}
\end{enumerate}	
{\rm Applying the functor $G_{\ta{(i_0)}}$ in the datum $\big(I,\gamma_{\ta(i_0)},\gamma_{(i_0)}\big)$ we obtain the following global action $\bt$ of $\Gamma_G^m$ on $A$:}

\begin{equation*}\bt_{(g,i,j)}=\begin{cases}
\text{\rm the identity map of } Ae_i,& \text{ if } (g,i,j)=(e,i,i),\\
\gamma_{\ta_{i}}\circ\sigma_{g}\circ\gamma\m_{\ta_{j}}\colon Ae_j\to Ae_i ,& \text{ otherwise } .
\end{cases}
\end{equation*}
\end{exe}

\vu

\section{Globalization}

\vu

In this section we study the globalization problem for a lifted partial groupoid action.
We start by recalling that a {\it globalization} of a partial action $\af=(A_g,\alpha_g)_{g\in \G}$ of a groupoid $\G$ on a ring $A$ is a pair $(\tilde{A},\tilde{\alpha})$, where $\tilde{A}$ is a ring and $\tilde{\alpha}=(\tilde{A}_g,\tilde{\alpha}_g)_{g\in \G}$ is a global action of $\G$ on $\tilde{A}$, that satisfies the following properties:
\begin{enumerate}[$\quad$(a)]
	\item  $A_y$ is an ideal of $\tilde{A}_{y}$, for all $y\in \G_0$,\vu
	\item  $A_g=A_{t(g)}\cap \tilde{\alpha}_g(A_{s(g)})$, \vu
	\item  $\tilde{\alpha}_g(a)=\alpha_g(a)$, for all $a\in A_{g^{-1}}$, \vu
	\item  $\tilde{A}_g=\sum_{t(h)=t(g)}\tilde{\alpha}_h(A_{s(h)})$, for all $g\in \G$.
\end{enumerate}

\vu

When $\alpha$ admits a globalization we say that $\alpha$ is {\it globalizable}. The globalization of $\alpha$ is unique, up to isomorphism; see \cite[Section 2]{BP} for details.

\begin{teo1}[\cite{BP}, Theorem 2.1]\label{teo:globaliza} Let $\alpha=\left(\{A_g\}_{g\in \G},\{\alpha_g\}_{g\in \G}\right)$ be a partial
	action of a groupoid $\G$ on a ring $A$ and suppose that $\A_y$ is a unital ring for each $y\in \G_0$. Then,
 $\alpha$ is globalizable  if and only if each ideal $A_g$, $g\in \G$, is a unital ring.
\end{teo1}

In what follows in this section, $\G$ denotes a connected groupoid, $x\in \G_0$ is a fixed object of $\G$ and  $\tau(x)=\{\tau_{y}:y\in \G_0\}$ is a transversal for $x$. We also fix a ring $A$, a datum $\gamma=(\I, \gamma_{\ta(x)},\gm_{(x)})\in \D_{\ta(x)}(\A)$ and $\bt=(B_g,\bt_g)_{g\in\G}$ the corresponding lifted partial groupoid action, that is, $\bt$ is the partial action of $\G$ on $A$ given by \eqref{afcon} and \eqref{dom}.

\vu

The next result relates the existence of globalization for $\beta$ with the existence of globalization for $\gamma_{(x)}$.

\begin{prop1}\label{prop:globa} The  partial groupoid action $\beta$ of $\G$ on $A$ is globalizable if and only if the partial group action $\gm_{(x)}$ of $\G(x)$ on $I_x$ is globalizable and the ideals $I_y,\,I_{\ta\m_{y}}$ are unital rings, for all $y\in \G_0$.
\end{prop1}
\begin{proof} Suppose that $\beta$ is globalizable. By Theorem \ref{teo:globaliza}, the ideal $B_g$ is a unital ring, for all $g\in \G$. Let $1_g$ be a central idempotent of $A$ such that $B_g=A1_g$, $g\in \G$. Then $B_h=A1_h=A1_x1_h=B_x1_h=I_x1_h$, for all $h\in \G(x)$. Thus, by Theorem 4.5 of \cite{DE}, the  partial group action $\gm_{(x)}$ of $\G(x)$ on $I_x$ is globalizable. Also, $I_y=B_y=A1_y$ and $I_{\ta\m_{y}}=B_{\ta\m_{y}}=A1_{\ta\m_{y}}$, for all $y\in \G_0$.
	
	Conversely, since $\gm_{(x)}$ is globalizable it follows from Theorem 4.5 of \cite{DE} that $I_h$ is a unital ring for all $h\in \G(x)$. Hence there are central idempotents $1_h$ of $I_x$ such that $I_h=I_x1_h$. Consider also central idempotents $1_y,\,1_{\ta\m_{y}}$ of $A$ such that $I_y=A1_y$ and $I_{\ta\m_{y}}=A1_{\ta\m_{y}}$ for all $y\in \G_0$. By \eqref{dom},   $B_{g}=\gamma_{\ta_{t(g)}}(\I_{\ta\m_{t(g)}}\cap \gamma_{g_x}(\I_{\ta\m_{s(g)}}\cap \I_{g\m_{x}}))$ for all $g\in \G$, $g\notin \G_0$. It is straightforward to check that
	$1_g=\gamma_{\ta_{t(g)}}(1_{\ta\m_{t(g)}}\gamma_{g_x}(1_{\ta\m_{s(g)}}1_{g\m_x}))$ is a central idempotent of $A$ that satisfies $A1_g=B_g$, for all $g\in \G$, $g\notin \G_0$. If $g=y\in \G_0$ then $B_g=B_y=I_y=A1_y$. Hence, Theorem \ref{teo:globaliza} implies that $\bt$ is globalizable.	
\end{proof}
\vu

Now we present explicitly the globalization of $\beta$ for a specific type of datum $\gamma$.
In order to do this we introduce the following definition. The datum $\gamma=(\I, \gamma_{\ta(x)},\gm_{(x)})\in \D_{\ta(x)}(\A)$ will be call {\it globalizable}  if
\begin{enumerate}[$\quad$(a)]
	\item $I_x$ is a unital ring and the partial group action $\gm_{(x)}$ of $\G(x)$ on $I_x$ admits a globalization $(J_x,\tilde{\gm}_{(x)})$,
	\item there exist a ring $B$, a family of ideals $\{J_y\}_{y\in\G_0}$ of $B$ and a family of isomorphisms $\{\tilde{\gamma}_{\ta_y}:J_x\to J_y\}_{y\in\G_0}$ of rings such that $I_y$ is an ideal of $J_y$ and $\tilde{\gamma}_{\ta_y} (a)=\gamma_{\ta_{y}}(a)$, for all $a\in I_{\ta\m_{y}}$ and $y\in \G_0$,
	\item  $I_{\ta\m_{y}}=I_x$ and $I_{\ta_{y}}=I_y$, for all $y\in \G_0$.
\end{enumerate}
We will denote: $\tilde{\gamma}=(J, \tilde{\gamma}_{\ta(x)},\tilde{\gm}_{(x)})$, where $J=\{J_y\}_{y\in \G_0}$ and $\tilde{\gamma}_{\ta(x)}=\{\tilde{\gamma}_{\ta_y}:J_x\to J_y\}_{y\in\G_0}$. The triple $\tilde{\gamma}$ will be call a {\it globalization} of $\gamma$ in $B$.

\vu

Suppose that $\gamma$ is globalizable. Then $\gamma_{\ta_y}: I_{\ta\m_y}=I_x\to I_{\ta\m_y}=I_y$ is an isomorphism and whence $I_x=I_{\ta\m_y}$ and $I_y=I_{\ta_y}$ are unital rings, for all $y\in\G_0$. Consequently, by Proposition \ref{prop:globa}, $\beta$ is globalizable. The next theorem describes the globalization of $\beta$ in the case that $\gamma$ is globalizable.

\begin{teo1}\label{teo:globaliza1} If $\tilde{\gamma}=(J, \tilde{\gamma}_{\ta(x)},\tilde{\gm}_{(x)})$ is a globalization of $\gamma$ in $B$ then the groupoid action $\tilde{\beta}$ of $\G$ on $B$ given by
	\[\bt=(\tilde{B}_g,\tilde{\bt}_g)_{g\in	\G},\quad \tilde{B}_g=J_{t(g)},\quad \tilde{\bt}_g=\tilde{\gamma}_{\ta_{t(g)}}\circ
	\tilde{\gamma}_{g_x}\circ \tilde{\gamma}\m_{\ta_{s(g)}}:J_{s(g)}\to J_{t(g)},\]
	is a globalization of $\beta$.
\end{teo1}
\begin{proof}
	Clearly $\tilde{\beta}$ is a global action of $\G$ on $B$ and $\tilde{\beta}_g(a)=\beta_g(a)$, for all $a\in B_{g^{-1}}$ and $g\in \G$. Notice that
	\begin{align*}
	B_{t(g)}\cap \tilde{\bt}_g(B_{s(g)})&= I_{t(g)}\cap\tilde{\bt}_g(I_{s(g)})\\
	&= I_{t(g)}\cap\tilde{\gamma}_{\ta_{t(g)}}\tilde{\gamma}_{g_x}\tilde{\gamma}\m_{\ta_{s(g)}}(I_{s(g)})\\
	&= I_{\ta_{t(g)}}\cap\tilde{\gamma}_{\ta_{t(g)}}\tilde{\gamma}_{g_x}\tilde{\gamma}\m_{\ta_{s(g)}}(I_{\ta_{s(g)}})\\
	&=\tilde{\gamma}_{\ta_{t(g)}}(I_x)\cap \tilde{\gamma}_{\ta_{t(g)}}\tilde{\gamma}_{g_x}(I_x)\\
	&=\tilde{\gamma}_{\ta_{t(g)}}(I_x\cap\tilde{\gamma}_{g_x}(I_x))\\
	&\overset{\mathclap{(\ast)}}{=}\tilde{\gamma}_{\ta_{t(g)}}(I_{g_x}),
	\end{align*}
	where $(\ast)$ holds because $(J_x,\tilde{\gm}_{(x)})$ is a globalization of the partial group action $\gm_{(x)}$ of $\G(x)$ on $I_x$. On the other hand,
	\begin{align*}
	B_g&=\gamma_{\ta_{t(g)}}(I_{\ta\m_{t(g)}}\cap\gamma_{g_x}(I_{\ta\m_{s(g)}}\cap I_{g\m_x}))\\
	&=\gamma_{\ta_{t(g)}}(I_{x}\cap\gamma_{g_x}(I_{x}\cap I_{g\m_x}))\\
	&=\gamma_{\ta_{t(g)}}(I_{g_x})\\
	&=\tilde{\gamma}_{\ta_{t(g)}}(I_{g_x}).
	\end{align*}
	Hence, $B_g=B_{t(g)}\cap \tilde{\bt}_g(B_{s(g)})$. Note also that, for each $g\in \G$, we have
	\begin{align*}
	\sum_{t(h)=t(g)}\tilde{\bt}_h(B_{s(h)})&=\sum_{t(h)=t(g)}\tilde{\gamma}_{\ta_{t(h)}}\tilde{\gamma}_{h_x}\tilde{\gamma}\m_{\ta_{s(h)}}(I_{s(h)})\\
	&=\tilde{\gamma}_{\ta_{t(g)}}\left(\sum_{t(h)=t(g)}\tilde{\gamma}_{h_x}(I_x)\right)\\
	&\overset{\mathclap{(\ast\ast)}}{=}\tilde{\gamma}_{\ta_{t(g)}}(J_x)\\
	&=J_{t(g)}\\
	&=\tilde{B}_g.
	\end{align*}
	In order to justify $(\ast\ast)$ note that the restriction of the map $\pi:\G\to \G(x)$ (see \eqref{epi}) to $\mathcal{T}_{t(g)}$ is surjective. In fact, given $l\in \G(x)$ we have that $\pi(\ta_{t(g)}l)=(\ta_{t(g)})_xl_x=xl=l$. Thus, $\sum_{t(h)=t(g)}\tilde{\gamma}_{h_x}(I_x)=\sum_{l\in \G(x)}\tilde{\gamma}_{l}(I_x)$. Since $(J_x,\tilde{\gm}_{(x)})$ is the globalization of $\gamma_{(x)}$ it follows that $J_x=\sum_{l\in \G(x)}\tilde{\gamma}_{l}(I_x)$.
\end{proof}

\vu

\begin{exe}\label{ex:datum-globa}
	{\rm Let $\G=\{g,h,l,m,l^{-1},m^{-1}\}\rightrightarrows \{x,y\}=\G_0$ be the groupoid given in Example \ref{hexd} and fix  $\ta(x)=\{\ta_x=x, \ta_y=l\}$. Let $A$ be a ring, $\sigma,\gamma$ ring automorphisms of $A$ with $\sigma\m=\sigma$ and $e\in A$ a central idempotent of $A$. Consider the datum $\big(I,\gamma_{\ta(x)},\gamma_{(x)}\big)$, where
	\begin{enumerate}[$\quad\circ $]
		\item  $I=\{I_x=Ae,\,I_y=A\gamma(e)\}$ is the set of the ideals of $A$,\vu
		\item $\gamma_{\ta(x)}=\{\gamma_{\ta_x}=\gamma_x=id_{I_x}, \gamma_{\ta_y}=\gamma_l=\gamma|_{I_x}:I_x\to I_y\}$ is the set of ring isomorphisms,\vu
		\item  $\gamma_{(x)}$ is the partial action of $\G(x)=\{x, g\}$ on $I_x$ given by $\gamma_x=id_{I_x}$ and \linebreak $\gamma_g=\gamma_{g\m}=\sigma|_{Ae\sigma(e)}:Ae\sigma(e)\to Ae\sigma(e)$.
	\end{enumerate}	
The lifted partial groupoid action corresponding to this datum is the partial action $\beta=(B_z,\bt_z\big)_{z\in\G}$ of $\G$ on $A$ given by
\begin{align*}
&\bt_x=id_{I_x},\quad \bt_y=id_{I_y},\quad \bt_g=\gamma_g,\quad \bt_l=\gamma_l,\quad\beta_m=\gamma\sigma|_{Ae\sigma(e)}:Ae\sigma(e)\to A\gamma(e\sigma(e)),\\
&\bt_h=\gamma\sigma \gamma\m|_{A\gamma(e\sigma(e))}:A\gamma(e\sigma(e))\to A\gamma(e\sigma(e)),\quad\beta_{l\m}=\gamma_l\m,\quad\beta_{m\m}=\sigma\gamma\m|_{A\gamma(e\sigma(e))}.
\end{align*}
The datum$\big(I,\gamma_{\ta(x)},\gamma_{(x)}\big)$ admits a globalization $\big(J,\tilde{\gamma}_{\ta(x)},\tilde{\gamma}_{(x)}\big)$ in $A$,  where
\begin{enumerate}[$\quad\circ $]
	\item  $J=\{J_x=J_y=A\}$ is the set of ideals of $A$,\vu
	\item  $\tilde{\gamma}_{\ta(x)}=\{\tilde{\gamma}_{\ta_x}=\tilde{\gamma}_x=id_{A}, \tilde{\gamma}_{\ta_y}=\gamma\}$ is the set of ring isomorphisms,\vu
	\item  $\tilde{\gamma}_{(x)}$ is the globalization of $\gamma_{(x)}$, i.~e. $\tilde{\gamma}_{(x)}$ is the global action of $\G(x)$ on $A$ given by $\tilde{\gamma}_{x}=id_{A}$ and $\tilde{\gamma}_{g}=\sigma$.
	\end{enumerate}		
By Theorem \ref{teo:globaliza1}, the global action $\tilde{\beta}$ of $\G$ on $A$ given by
\[
\tilde{\beta}_{x}=\tilde{\beta}_{y}=id_A,\,\,\, \tilde{\beta}_{g}=\sigma,\,\,\, \tilde{\beta}_{h}=\gamma\sigma\gamma\m,\,\,\,\tilde{\beta}_{m}=\gamma\sigma,\,\,\, \tilde{\beta}_{m\m}=\sigma\gamma\m,\,\,\,\tilde{\beta}_{l}=\gamma,\,\,\,\tilde{\beta}_{l\m}=\gamma\m,
\]
is the globalization of the partial action $\beta$.}
\end{exe}

\section{The partial skew groupoid ring}

\vu

In this section the partial action $\af=(A_g,\af_g)_{g\in \G}$ of $\G$ on the ring $A$ is assumed to be unital, that is, each $A_g$ is unital with identity element denoted by $1_g$ which is a central idempotent of $A$ and $A_g=A1_g$. Also, $x\in \G_0$ is a fixed object of $\G$ and $\ta(x)=\{\ta_y\ |\ y\in \G_0\}$ is a transversal for $x$.

\vu

Suppose that $\G_0$ is finite and $\af$ is $\ta(x)$-global. Our main goal in this section is to prove that  there exists  a unital ring $C$ and a partial  group action $\theta$   of $\G(x)$ on $C$ such that $A\star_\af\G$ is isomorphic to $C\star_ \theta\G(x)$.

\vu

We start by observing that $\G$ and $\G_0^2\times \G(x)$ are isomorphic as groupoids, where  the groupoid structure in $\G_0^2=\G_0\times\G_0$ is  $s(y,z)=y$ and $t(y,z)=z$, for all $(y,z)\in\G_0^2$.
By Lemma \ref{lema31}  the map $f:\G\to \G_0^2\times\G(x)$ defined by $g\mapsto ((s(g),t(g)), g_x)$, with $g_x=\ta^{-1}_{t(g)}g\ta_{s(g)}$, for all $g\in \G,$ is a groupoid epimorphis.
Suppose that $f(g)$ is an identity of $\G_0^2\times\G(x)$. Then $f(g)=((y,y),x)$ for some $y\in\G_0$. Hence $s(g)=t(g)=y$ and $x=g_x=\ta_y^{-1}g\ta_y$ which implies that $g=\ta_y\ta_y^{-1}=x$. This ensures that $f$ is injective.  \vu

\begin{lem1}\label{xx}
	If $\af$ is $\tau(x)$-global then $\af^\ast=(A^\ast_u,\af^\ast_u)_{u\in\G_0^2}$, given by $A^\ast_u=A_{t(u)}$ and \linebreak $\aft_u=\af_{\ta_{t(u)}}\circ\af_{\ta\m_{s(u)}}$, is a global action of $\G_0^2$ on $A$.
\end{lem1}

\begin{proof}
	Indeed, for any identity $e=(y,y)$ of $\G_0^2$ we have that $A^\ast_e=A_y$ and $\aft_e=\af_{\ta_{y}}\circ\af_{\ta\m_{y}}$ is the identity map of $A_{\ta_y}=A_y.$ Moreover, if  $u=(y,z)$ and $v=(r,w)$ are elements in $\G_0^2$ such that the product $uv$  exists, then $uv=(r,z)$ and  $ \aft_u\circ \aft_v=\af_{\ta_z}\circ\af_{\ta\m_y}\circ\af_{\ta_y}\circ\af_{\ta\m_r}=\af_{\ta_z}\circ\af_{\ta\m_r}= \aft_{uv}$.
\end{proof}

Recalling from \cite[Section 3]{BP}, that the partial skew groupoid ring $A\star_\af\G$ is
defined as the direct sum
$$A\star_\alpha \G=\bigoplus_{g\in \G}A_g\delta_g$$ (where the $\delta_g$'s are placeholder symbols) with the usual addition and multiplication
induced by the rule
\[
(a\delta_g)(b\delta_h)=
\begin{cases}
a\alpha_g(b1_{g^{-1}})\delta_{gh} &\text{if $s(g)=t(h)$}\\
0  &\text{otherwise},
\end{cases}
\]
for all $g, h\in \G$, $a\in A_g$ and $b\in A_h$. Endowed with such a multiplication $A\star_\af\G$ has a structure of an associative ring. If in addition the set $\G_0$ of the objects of $\G$ is finite then $A\star_\af\G$ is also unital with identity element $1_{A\star_\af\G}=\sum_{y\in\G_0}1_y\delta_y$.

\vu

Thanks to Lemma \ref{xx} we can  consider the corresponding (global) skew groupoid ring $C=A\star_ {\aft}\G_0^2$.
In the sequel we will also see, under the same conditions listed in Lemma \ref{xx}, that the group $\G(x)$ acts on $C$ via an appropriate partial action $ \theta$, which will allows us to consider the corresponding partial skew group ring $C\star_ \theta\G(x)$.

\vu

In order to construct  $ \theta$ we start by looking for a family of ideals in $C$. For this purpose set
\begin{align}\label{Czh}
C_{z,h}=\af_{\ta_z}(A_h),\,\ \text{for all}\,\  z\in\G_0\,\,\text{and}\,\ h\in\G(x),
\end{align}
and take
\begin{align}\label{Ch}
C_h=\oplus_{u\in\G_0^2}C_{t(u),h}\delta_u.
\end{align}

\begin{lem1}\label{63} If $\af$ is $\tau(x)$-global then the following statements are true:
	\begin{enumerate}[\rm (i)]
		\item $C_x=C$, \vu
		\item If $\G_0$ is finite then $C_h$ is a unital ideal of $C$, for all $h\in\G(x)$.
	\end{enumerate}

\end{lem1}
\begin{proof}
	Firstly,
	\begin{align*}
	C &=\oplus_{u\in \G_0^2}A_{{t(u)}}\delta_u\\
	  &= \oplus_{u\in \G_0^2}A_{\ta_{t(u)}}\delta_u \\
      &= \oplus_{u\in \G_0^2}\af_{\ta_{t(u)}}(A_{\ta\m_{t(u)}})\delta_u \\
      &= \oplus_{u\in \G_0^2}\af_{\ta_{t(u)}}(A_x)\delta_u \\
      &=C_x.
      \end{align*}
Since $\G_0$ is finite, $1'_h=\sum_{z\in\G_0}\af_{\ta_z}(1_h)\delta_{(z,z)}$ is the identity element of $C_h$, for all $h\in\G(x)$. Indeed, for every $a=\af_{\ta_w}(a_h)\delta_{(y,w)}\in C_{w,h}\delta_u$, with $a_h\in A_h,$ one has that
\begin{align*}
a1'_h &=\sum_{z\in\G_0}\af_{\ta_w}(a_h)\delta_{(y,w)}\af_{\ta_z}(1_h)\delta_{(z,z)}\\
&=\af_{\ta_w}(a_h)\delta_{(y,w)}\af_{\ta_y}(1_h)\delta_{(y,y)}\\
&=\af_{\ta_w}(a_h) \aft_{(y,w)}(\af_{\ta_y}(1_h))\delta_{(y,w)}\\
&\overset{\mathclap{(\ast)}}{=}\af_{\ta_w}(a_h)\af_{\ta_w}\af_{\ta\m_y}\af_{\ta_y}(1_h)\delta_{(y,w)}\\
&=\af_{\ta_w}(a_h)\af_{\ta_w}(1_h)\delta_{(y,w)}\\
&=\af_{\ta_w}(a_h)\delta_{(y,w)}\\
&=a,
\end{align*}
where $(\ast)$ follows from Lemma \ref{xx}. By similar arguments one gets $1'_ha=a$. Hence $1'_h$ is a central idempotent of $C$ and it is easy to see that $C_h=1'_hC$.
\end{proof}

\vu

For each $(z,h)\in \G_0\times \G(x)$ define $ \theta_{z,h}:C_{z,h\m}\to C_{z,h}$ by $\af_{\ta_z}(a)\mapsto \af_{\ta_z}(\af_h(a))$, for all $a\in A_{h\m}$. Clearly $ \theta_{z,h}$ is a bijection. Moreover, these maps induce the bijection $ \theta_h:C_{h\m}\to C_h$, given by  $ \theta_h(\af_{\ta_{t(u)}}(a)\delta_u)= \theta_{t(u),h}(a)\delta_u$, for all $a\in A_{h\m}$ and $u\in \G_0^2$.

\vu

In all what follows in this section we will assume that $\G_0$ is finite.

\begin{lem1}\label{64}
The pair	$\theta=(C_h, \theta_h)_{h\in \G(x)}$ is a partial action of $\G(x)$ on $C$.
\end{lem1}	

\begin{proof}
We will proceed by steps.

\vu

\noindent \emph{Step 1: $ \theta_x$ is the identity map of $C_x=C$.}

\vu

It  follows straightforward from the definition of $ \theta$.

\vd

\noindent \emph{Step 2: $ \theta_h$ is multiplicative, for all $h\in\G(x)$.}

\vu

For all $a,b\in A_{h\m}$,
\begin{align*}
\hspace*{-1.5cm} \theta_h((\af_{\ta_z}(a)\delta_{(y,z)})(\af_{\ta_y}(b)\delta_{(w,y)}))&= \theta_h(\af_{\ta_z}(a) \aft_{(y,z)}(\af_{\ta_y}(b))\delta_{(w,z)})\\
&= \theta_h(\af_{\ta_z}(a)\af_{\ta_z}(\af_{\ta\m_y}(\af_{\ta_y}(b)))\delta_{(w,z)})\\
&= \theta_h(\af_{\ta_z}(ab)\delta_{(w,z)})\\
&=\af_{\ta_z}(\af_h(ab))\delta_{(w,z)}\\
&=\af_{\ta_z}(\af_h(a))\af_{\ta_z}(\af_h(b))\delta_{(w,z)}.
\end{align*}
On the other hand,
\begin{align*}
\hspace*{-1.2cm} \theta_h(\af_{\ta_z}(a)\delta_{(y,z)}) \theta_h(\af_{\ta_y}(b)\delta_{(w,y)})&=\af_{\ta_z}(\af_h(a))\delta_{(y,z)}\af_{\ta_y}(\af_h(b))\delta_{(w,y)}\\
&=\af_{\ta_z}(\af_h(a))\af_{\ta_z}(\af_{\ta\m_y}(\af_{\ta_y}(\af_h(b))))\delta_{(w,z)}\\
&=\af_{\ta_z}(\af_h(a))\af_{\ta_z}(\af_h(b))\delta_{(w,z)}.
\end{align*}

\vu

\noindent  \emph{Step 3: $ \theta_{l\m}(C_l\cap C_{h\m})\subset C_{(hl)\m}$, for all $h,l\in \G(x)$.}
\begin{align*}
\hspace*{-1cm} \theta_{l\m}(C_l\cap C_{h\m})&=\oplus_{u\in \G_0^2} \theta_{l\m}(\af_{\ta_{t(u)}}(A_l\cap A_{h\m})\delta_u)\\
&=\oplus_{u\in \G_0^2}\af_{\ta_{t(u)}}(\af_{l\m}(A_l\cap A_{h\m}))\delta_u\\
&\subset \oplus_{u\in \G_0^2}\af_{\ta_{t(u)}}(A_{l\m}\cap A_{l\m h\m})\delta_u\quad (\text{for}\,\,\af_{(x)}\,\,\text{is a partial group action})\\
&=C_{(hl)\m}.
\end{align*}
\noindent  \emph{Step 4: $ \theta_{l} ( \theta_h (c))= \theta_{hl}(c)$, for all $c=\af_{\ta_z}(\af_{l\m}(a))\delta_{(y,z)}\in  \theta_{l\m}(C_l\cap C_{h\m})$.}
\begin{align*}
 \theta_h( \theta_l(c))&= \theta_h( \theta_l(\af_{\ta_z}(\af_{l\m}(a))\delta_{(y,z)}))\\
&= \theta_h(\af_{\ta_z}(a)\delta_{(y,z)})\\
&=\af_{\ta_z}(\af_h(a))\delta_{(y,z)}.
\end{align*}
On the other hand, since $\af_{hl}=\af_h\af_l$ in $\af_{l\m}(A_l\cap A_{h\m})$ we have
\[ \theta_{hl}(c)=\af_{\ta_z}(\af_{hl}(\af_{l\m}(a)))\delta_{(y,z)}=\af_{\ta_z}(\af_h(a))\delta_{(y,z)}.\]
Hence $ \theta_{hl}= \theta_h \theta_l$ in $ \theta_{l\m}(C_l\cap C_{h\m})$.
\end{proof}

Now we will prove the main result of this section.

\begin{teo1} \label{la grande finale}
	If $\af$ is $\tau(x)$-global then
	$$\varphi:A\star_\af\G\to (A\star_{\aft}\G_0^2)\star_ \theta\G(x),\quad a\delta_g\mapsto a\delta_{(s(g),t(g))}\delta_{g_x}$$ is a ring isomorphism.	
\end{teo1}

\begin{proof}

	We proceed again by steps.
	
	\vu
	
\noindent	\emph{Step 1: $\varphi$ is well defined.}

\vu
	
First of all we have by Lemma \ref{xx} that $A_g\subseteq A_{{t(g)}}=A^\ast_{(s(g),t(g))}$, for all $g\in\G$. Hence,  we only need to show that $a\delta_{(s(g),t(g))}\in C_{g_x}$, for all $a\in A_g$.
Notice that
	\begin{align*}
	\af_{\ta_{t(g)}}(A_{g_x})&=\af_{\ta_{t(g)}}(A_{\ta\m_{t(g)}g\ta_{s(g)}})\\
	&=\af_{\ta_{t(g)}}(A_{\ta\m_{t(g)}g\ta_{s(g)}}\cap A_x)\\
	&\overset{\mathclap{\eqref{bapa3}}}{=}A_{g\ta_{s(g)}}\cap A_{\ta_{t(g)}x}\\
	&=A_{g\ta_{s(g)}}\cap A_{\ta_{t(g)}}\\
    &\overset{\mathclap{\eqref{cond1}}}{=}A_{g\ta_{s(g)}}\cap A_{t(g)}\qquad (A_{g\ta_{s(g)}}\subseteq A_{t(g\ta_{s(g)})}=A_{t(g)} )\\
	&=A_{g\ta_{s(g)}}.
	\end{align*}
Since $A_{g\m}=A_{g\m}\cap A_{s(g)}\overset{\eqref{cond1}}{=} A_{g\m}\cap A_{\ta_{s(g)}}$ we have $A_g=\af_g(A_{g\m}\cap A_{\ta_{s(g)}})=A_g\cap A_{g\ta_{s(g)}}$. Hence $a\in A_g\subseteq A_{g\ta_{s(g)}}=\af_{\ta_{t(g)}}(A_{g_x})$ which implies $a\delta_{(s(g),t(g))}\in C_{g_x}$ by \eqref{Czh} and \eqref{Ch}.

\vd

\noindent \emph{Step 2: $\varphi$ is a ring homomorphism.}

\vu	

For all $g,h\in\G$ such that $s(g)=t(h)$, $a\in A_g$ and $b\in A_h$  we have
\begin{align*}
\varphi((a\delta_g)(b\delta_h))&=\varphi(a\af_g(b1_{g\m})\delta_{gh})\\
&=a\af_g(b1_{g\m})\delta_{(s(gh),t(gh))}\delta_{(gh)_x}\\
&=a\af_g(b1_{g\m})\delta_{(s(h),t(g))}\delta_{g_xh_x}.
\end{align*}
On the other hand
\begin{align*}
\varphi(a\delta_g)\varphi(b\delta_h)&=(a\delta_{(s(g),t(g))}\delta_{g_x})(b\delta_{(s(h),t(h))}\delta_{h_x})\\
&=(a\delta_{(s(g),t(g))})(\gamma_{g_x}(b\delta_{(s(h),t(h))}1'_{g\m_x}))\delta_{g_xh_x}.
\end{align*}
Since $b\in A_h\overset{\eqref{cond1}}{\subseteq} A_{h\ta_{s(h)}}=\af_{\ta_{t(h)}}(A_{h_x})$ there is $b'\in A_{h_x}$ such that $b=\af_{\ta_{t(h)}}(b')$ and
\begin{align*}
b\delta_{(s(h),t(h))}1'_{g\m_x}&=\af_{\ta_{t(h)}}(b')\delta_{(s(h),t(h))}\sum_{g\in\G_0}\af_{\ta_z}(1_{g\m_x})\delta_{(z,z)}\\
&=\af_{\ta_{t(h)}}(b')\delta_{(s(h),t(h))}\af_{\ta_{s(h)}}(1_{g\m_x})\delta_{(s(h),s(h))}\\
&=\af_{\ta_{t(h)}}(b')\af_{\ta_{t(h)}}(1_{g\m_x})\delta_{(s(h),t(h))}\\
&=\af_{\ta_{t(h)}}(b'1_{g\m_x})\delta_{(s(h),t(h))}.
\end{align*}
Hence,
\begin{align*}
\varphi(a\delta_g)\varphi(b\delta_h)&=(a\delta_{(s(g),t(g))})(\gamma_{g_x}(b\delta_{(s(h),t(h))}1'_{g\m_x}))\delta_{g_xh_x}\\
&=(a\delta_{(s(g),t(g))})(\gamma_{g_x}(\af_{\ta_{t(h)}}(b'1_{g\m_x})\delta_{(s(h),t(h))}))\delta_{g_xh_x}\\
&=(a\delta_{(s(g),t(g))})(\af_{\ta_{t(h)}}(\af_{g_x}(b'1_{g\m_x}))\delta_{(s(h),t(h))})\delta_{g_xh_x}\\
&=a\af_{\ta_{t(g)}}\af_{\ta\m_{s(g)}}\af_{\ta_{t(h)}}\af_{g_x}(b'1_{g\m_x})\delta_{(s(h),t(g))}\delta_{g_xh_x}\\
&=a\af_{\ta_{t(g)}}\af_{g_x}(b'1_{g\m_x})\delta_{(s(h),t(g))}\delta_{g_xh_x}\\
&\overset{\mathclap{\eqref{cond1}}}{=}a\af_{\ta_{t(g)}}(\af_{g_x}(b'1_{g\m_x})1_{\ta\m_{t(g)}})\delta_{(s(h),t(g))}\delta_{g_xh_x}\\
&\overset{\mathclap{\eqref{par2}}}{=}a\af_{\ta_{t(g)}g_x}(b'1_{(\ta_{t(g)}g_x)\m})1_{\ta_{t(g)}}\delta_{(s(h),t(g))}\delta_{g_xh_x}\\
&=a\af_{g\ta_{s(g)}}(b'1_{(g\ta_{s(g)})\m})1_{\ta_{t(g)}}\delta_{(s(h),t(g))}\delta_{g_xh_x}\\
&\overset{\mathclap{\eqref{cond1}}}{=}a\af_{g\ta_{s(g)}}(b'1_{(g\ta_{s(g)})\m})1_g\delta_{(s(h),t(g))}\delta_{g_xh_x}\\
&\overset{\mathclap{\eqref{par2}}}{=}a\af_g(\af_{\ta_{s(g)}}(b'1_{\ta\m_{s(g)}})1_{g\m})\delta_{(s(h),t(g))}\delta_{g_xh_x}\\
&\overset{\mathclap{\eqref{cond1}}}{=}a\af_g(\af_{\ta_{t(h)}}(b')1_{g\m})\delta_{(s(h),t(g))}\delta_{g_xh_x}\\
&=a\af_g(b1_{g\m})\delta_{(s(h),t(g))}\delta_{g_xh_x}\\
&=\varphi((a\delta_g)(b\delta_h)).
\end{align*}

\vd

\noindent \emph{Step 3: $\varphi$ is injective.}

\vu

Let $v=\sum_{g\in\G}a_g\delta_g\in A\star_\af\G$ such that $\varphi(v)=0$. Then
$$0=\sum_{g\in\G}a_g\delta_{(s(g),t(g))}\delta_{g_x}=\sum_{h\in\G(x)}\sum_{\pi(g)=h}a_g\delta_{(s(g),t(g))}\delta_h$$ which implies that $\sum_{\pi(g)=h}a_g\delta_{(s(g),t(g))}=0$ ($\ast$) for all $h\in\G(x)$ (notice that $C\star_ \theta\G(x)$ is a direct sum).	Finally, it is straightforward to check that for any $h\in\G(x)$ and any $g,g'\in\pi\m(h)$, $(s(g),t(g))=(s(g'),t(g'))$ if and only if $g=g'$. Therefore ($\ast$) holds if and only if $a_g=0$, for all $g\in \G$, and so $v=0$.

\vd

\noindent \emph{Step 4: $\varphi$ is surjective.}

\vu

It is enough to check that given any element of the type $\af_{\ta_z}(a)\delta_{(y,z)}\delta_h$, with $h\in\G(x)$ and $a\in A_h$, there exists an element $w\in A\star_\af\G$ such that $\varphi(w)=\af_{\ta_z}(a)\delta_{(y,z)}\delta_h$. To do that observe firstly that
\begin{align*}
\af_{\ta_z}(a)\in\af_{\ta_z}(A_h)&=\af_{\ta_z}(A_h\cap A_x)\\
&\overset{\mathclap{\eqref{cond1}}}{=}\af_{\ta_z}(A_h\cap A_{\ta\m_z})\\
&\overset{\mathclap{\eqref{bapa3}}}{=}A_{\ta_zh}\cap A_{\ta_z}\\
&\overset{\mathclap{\eqref{cond1}}}{=}A_{\ta_zh}\cap A_z\\
&=A_{\ta_zh}\qquad  (\text{for}\,\,\ A_{\ta_zh}\subseteq A_{t(\ta_zh)}=A_{t(\ta_z)}=A_z).
\end{align*}
Therefore, for  $g=\ta_zh\ta\m_y$ we have $t(g)=t(\ta_z)=z$, $s(g)=s(\ta\m_y)=y$ and  $$\af_{\ta_z}(a)\in A_{\ta_zh}= A_{\ta_{t(g)}h}\overset{(\ast)}{=}\af_{\ta_{t(g)}}(A_h)\overset{(\ast\ast)}{\subseteq}\af_{\ta_{t(g)}}(A_{h\ta\m_{s(g)}})\overset{(\ast)}{=} A_{\ta_{t(g)}h\ta\m_{s(g)}}=A_{\ta_zh\ta\m_y}=A_g,$$ where $(\ast)$ is ensured by:
$$\af_{\ta_{t(g)}}(A_h)=\af_{\ta_{t(g)}}(A_h\cap A_x)\overset{\eqref{cond1}}{=}\af_{\ta_{t(g)}}(A_h\cap A_{\ta\m_{s(g)}})=A_{\ta_{t(g)}h},$$ and $(\ast\ast)$ by:
$$A_h=\af_h(A_{h\m}\cap A_x)\overset{\eqref{cond1}}{=}\af_h(A_{h\m}\cap A_{\ta\m_{s(g)}})=A_h\cap A_{h\ta\m_{s(g)}}\subseteq A_{h\ta\m_{s(g)}}.$$

\nod Now, taking $w=\af_{\ta_z}(a)\delta_g$ we are done.	
\end{proof}

The next proposition characterizes when a lifted partial action is $\tau(x)$-global.

\begin{prop1}\label{58}
	A lifted partial action of $\G$ on $A$ coming from  a  datum $\big(I,\gamma_{\ta(x)},\gamma_{(x)}\big)$ in $\D_{\ta(x)}(A)$ is $\ta(x)$-global if and only if $I_{\ta\m_y}=I_x$ and $I_{\ta_y}=I_y$ for all $y\in\G_0$.
\end{prop1}

\begin{proof}
	Along the proof we will freely use the formula \eqref{dom} without any mention to it.
	Let $\bt=\big(B_g, \bt_g\big)$ be the lifted partial action coming from $\big(I,\gamma_{\ta(x)},\gamma_{(x)}\big)$. If $ \bt$ satisfies \eqref{cond1} then $B_{\ta\m_y}=B_x$, for  all $y\in\G_0$. Thus,
	\begin{align*}
	B_{\ta\m_y}&=\gamma_{\ta_{t(\ta\m_y)}}\big(I_{\ta\m_{t(\ta\m_y)}}\cap\gamma_{(\ta\m_y)_x}(I_{\ta\m_{s(\ta\m_y)}}\cap I_{(\ta\m_y)\m_x})\big)\\
	&=\gamma_x\big(I_x\cap\gamma_x(I_{\ta\m_y}\cap I_x)\big)=\gamma_x(I_{\ta\m_y})=I_{\ta\m_y},
	\end{align*}
	for  all $y\in\G_0$ with $y\neq x$. Hence, $I_{\ta\m_y}=B_{\ta\m_y}=B_x=I_x$ for  all $y\in\G_0$.
	It is analogous to verify that $I_{\ta_y}=I_y$.
	\vu
	
	Conversely, if $I_{\ta\m_y}=I_x$ and $I_{\ta_y}=I_y$ for all $y\in\G_0$, then $B_g=\gamma_{\ta_{t(g)}} (I_{g_x})$ for all $g\in \G$, $g\notin \G_0$. Particularly, $B_{\ta\m_y}=I_x=B_x$ and $B_{\ta_y}=I_{\ta_y}=I_y=B_y$ and whence $ \bt$ satisfies \eqref{cond1}.
\end{proof}

\vu

Now we introduce a definition that will be useful in the rest of the paper. A lifted partial action $\bt=(B_g,\bt_g)_{g\in\G}$ of $\G$ on $A$ coming from a datum $\big(I,\gamma_{\ta(x)},\gamma_{(x)}\big)\in \D_{\ta(x)}(\A)$ will be called {\it $\gamma$-unital} if $I_x$ is a unital ring and $\gamma_{(x)}$ is a unital partial group action.

\begin{remark}\label{ob-partial-unital} {\rm Let $\beta$ be a lifted partial action of $\G$ on $A$ which is $\ta(x)$-global and $\gamma$-unital. Then $\beta$ is unital, i.~e. $B_g$ is a unital ring for all $g\in \G$. Indeed, suppose that $I_x=A1_x$ and $I_h=A1_h$, with $1_x$ and $1_h$ central idempotents of $A$, for all $h\in \G(x)$.  Using that $\beta$ is $\ta$-global we obtain that  $B_g=\gamma_{\ta_{t(g)}}(I_{g_x})$ and whence $B_g=A\gamma_{\ta_{t(g)}}(1_{g_x})$, for all $g\in \G$. Thus $B_g=A1_g$, where $1_g=\gamma_{\ta_{t(g)}}(1_{g_x})$ is a central idempotent of $A$, for all $g\in \G$.}
\end{remark}

\vu

Throughout the rest of this paper we will assume that $\G$ is finite,  $\bt=(B_g,\bt_g)_{g\in\G}$  will denote a fixed partial action of $\G$ on $A$ lifted from a given  datum $\big(I,\gamma_{\ta(x)},\gamma_{(x)}\big)$ in $\D_{\ta(x)}(\A)$. We will also assume that $\beta$ is $\tau(x)$-global,  $\gamma$-unital and $A=\oplus_{y\in\G_0}I_{y}=\oplus_{y\in\G_0}B_{y}$.
As in Remark \ref{ob-partial-unital}, the unit of $B_g$ will be denoted by $1_g$. Precisely, $B_x=I_x=A1_x$, $B_h=I_h=A1_h$ and $B_g=A1_g$ where $1_g=\gamma_{\ta_{t(g)}}(1_{g_x})$, for all $h\in \G(x)$ and $g\in \G$.

\vu

\section{The subring of invariants and the trace map}

According to \cite{BP} an element $a\in A$ is called \emph{invariant} by $\bt$ if $\bt_g(a1_{g\m})=a1_g$, for all $g\in\G$. We will denote by $\A^\bt$ the set of all invariants of $A$. Clearly, $A^\bt$ is a subring of $A$.

\begin{prop1}\label{inv}
Let $b=\somay b_y\in A$.  Then $b\in A^\bt$ if and only if $b_x\in B_x^{\gamma_{(x)}}$ and $b_y=\gamma_{\ta_y}(b_x)$, for all $y\in\G_0$.
\end{prop1}

\begin{proof}
If $b\in A^\bt$ then $b_x1_h=b1_h=\bt_h(b1_{h\m})=\gamma_h(b_x1_{h\m})$ for all $h\in \G(x)$. Hence $b_x\in B_x^{\gamma_{(x)}}$.
Also, $\gamma_{\ta_y}(b_x)=\gamma_{\ta_y}\gamma_x\gamma\m_x(b_x1_x)=\bt_{\ta_y}(b1_{\ta\m_y})=b1_{\ta_y}=b_y$, for all $y\in \G_0$.

Conversely, let $b=\somay\gamma_{\ta_y}(b_x)$, for some $b_x\in B_x^{\gamma_{(x)}}$. Then, for all $g\in\G$,
\begin{align*}
\bt_g(b1_{g\m})&=\gamma_{\ta_{t(g)}}\gamma_{g_x}\gamma\m_{\ta_{s(g)}}(b1_{g\m})\\
&=\gamma_{\ta_{t(g)}}\gamma_{g_x}\gamma\m_{\ta_{s(g)}}(\gamma_{\ta_{s(g)}}(b_x)\gamma_{\ta_{s(g)}}(1_{g\m_x}))\quad (\text{by assumption})\\
&=\gamma_{\ta_{t(g)}}\gamma_{g_x}(b_x1_{g\m_x})\\
&=\gamma_{\ta_{t(g)}}(b_x1_{g_x})\\
&=b1_g.
\end{align*}\end{proof}
The \emph{$\bt$-trace map} is defined as the map $t_\bt:A\to A$ given by $t_\bt(a)=\somag \bt_g(a1_{g\m})$, for all $a\in A$. Similarly, the $\bt_{(x)}$-trace map corresponding to the partial  action $\bt_{(x)}$ of $\G(x)$ over $B_x$ is defined as the map $t_{\bt_{(x)}}:B_x\to B_x$ given by $t_{\bt_{(x)}}(b)=\sum_{h\in\G(x)}\bt_h(b1_{h\m})$, for all $b\in B_x$. Notice that, by Remark \ref{r32},  $t_{\bt_{(x)}}=t_{\gamma_{(x)}}$ .

It follows from \cite[Lemma 4.2]{BP} (resp., \cite[Lemma 2.1]{DFP}) that $t_\bt(A)\subseteq A^\bt$ (resp., $t_{\bt_{(x)}}(B_x)\subseteq B_x^{\bt_{(x)}}$)

\begin{prop1}\label{p52} Let $b_x\in B_x$ and $b_z=\gamma_{\ta_z}(b_x)\in B_z$. Then \[t_{\beta}(b_z)=\sum_{y\in \G_0} \gamma_{\tau_y}(t_{\beta_{(x)}}(b_x)).\]
\end{prop1}
\begin{proof}
Indeed,
\begin{align*}
t_{\beta}(b_z)&=\sum_{g\in \G}\beta_g(b_z1_{g\m})\\
              &\overset{\mathclap{\eqref{afcon}}}{=}\sum_{g\in \G}\gamma_{\tau_{t(g)}}\gamma_{g_x}\gamma\m_{\tau_{s(g)}}(\gamma_{\tau_{z}}(b_x)\gamma_{\tau_{s(g)}}(1_{g\m_x}))\\
              &=\sum_{g\in \mathcal{S}_z}\gamma_{\tau_{t(g)}}\gamma_{g_x}(b_x1_{g\m_x})\\
              &=\sum_{l\in \mathcal{S}_x}\gamma_{\tau_{t(l)}}\gamma_{l_x}(b_x1_{l\m_x}) \qquad (l=g\ta_z)\\
              &=\sum_{y\in \G_0}\sum_{l\in \G(x,y)}\gamma_{\tau_{y}}\gamma_{l_x}(b_x1_{l\m_x}) \\
              &=\sum_{y\in \G_0}\sum_{h\in \G(x)}\gamma_{\tau_{y}}\gamma_{(\ta_yh)_x}(b_x1_{(\ta_yh)\m_x}) \qquad (\ta_yh=l) \\
              &=\sum_{y\in \G_0}\sum_{h\in \G(x)}\gamma_{\tau_{y}}\gamma_{h}(b_x1_{h\m}) \\
              &=\sum_{y\in \G_0}\gamma_{\tau_{y}}(\sum_{h\in \G(x)}\gamma_{h}(b_x1_{h\m}))\\
              &=\sum_{y\in \G_0}\gamma_{\tau_{y}}(t_{\beta_{(x)}}(b_x)).
\end{align*}	
\end{proof}

\begin{cor1}\label{53}
Suppose that $a=\sum_{z\in\G_0}b_z\in A$, with $b_z=\gamma_{\ta_z}(c_z)$ for some %unique 
$c_z\in B_x$, for each $z\in\G_0$. Then, $t_\bt(a)=\sum_{y\in\G_0}\gamma_{\ta_y}(t_{\bt_{(x)}}(c))$, where
$ c=\somaz c_z$.
\end{cor1}

\begin{proof}
Indeed,
\begin{align*}
t_\bt(a)&=\somaz t_\bt(b_z)\\
&\overset{\text{Prop.}\ref{p52}}{=}\somaz\somay\gamma_{\ta_y}(t_{\bt_{(x)}}(c_z))\\
&=\somay\gamma_{\ta_y}\left(t_{\bt_{(x)}}\left(\somaz c_z\right)\right)
\\&=\somay\gamma_{\ta_y}(t_{\bt_{(x)}}(c)).
\end{align*}
\end{proof}

\begin{teo1} \label{t54}
$t_\bt(A)=A^\bt$ if and only if $t_{\bt_{(x)}}(B_x)=B_x^{\bt_{(x)}}$.
\end{teo1}

\begin{proof}
Assume that $t_\bt$ is onto. Take $b_x\in B_x^{\bt_{(x)}}$ and set $b=\somay\gamma_{\ta_y}(b_x)$ which lies in $A^\bt$ by Proposition \ref{inv}.
Hence $b=t_\bt(a)$ for some $a\in A$. Since $A=\oplus_{z\in\G_0}B_z$ and $\gamma_{\ta_z}$  is a ring isomorphism  from $B_x$ to $B_z$ for each $z\in\G_0$, it follows
that $a=\somaz \gamma_{\ta_z}(c_z)$, with $c_z\in B_x$. Therefore $\somay\gamma_{\ta_y}(b_x)=b=t_\bt(a)\overset{\text{Cor.}\ref{53}}{=}\somay\gamma_{\ta_y}(t_{\bt_{(x)}}(c))$ with $c=\somaz c_z$, which implies
that $b_x=t_{\bt_{(x)}}(c)\in t_{\bt_{(x)}}(B_x)$.

\vu

Conversely, it is enough to observe that, by Propositions \ref{inv} and \ref{p52}, and the assumption on $\bt_{(x)}$, each element in $A^\bt$ is of the form $b=\somay\gamma_{\ta_y}(b_x)=t_\bt(a_x),$ with $a_x\in B_x$ such that $t_{\bt_{(x)}} (a_x)=b_x$.
\end{proof}

\section{The related Morita theory}

We start by recalling the definition of a Morita context. Given two unital rings $R$ and $S$,
bimodules $_RU_S$ and $_SV_R$, and bimodule maps $\mu: U\otimes_S V\to R$ and $\nu:V\otimes_R U\to S$,
the sixtuple $(R,S,U,V,\mu,\nu)$ is called a \emph{Morita context} (\emph{associated with $_RU$}) if the set
$$\left(
\begin{array}{cc}
R & U \\
V & S \\
\end{array}
\right)=\left\{\left(
\begin{array}{cc}
r & u \\
v & s \\
\end{array}
\right)\ \big|\ r\in R, s\in S, u\in U, v\in V\right\}
$$ with the usual addition and multiplication given by the rule
$$\left(
\begin{array}{cc}
r & u \\
v & s \\
\end{array}
\right)\left(
\begin{array}{cc}
r' & u' \\
v' & s' \\
\end{array}
\right)=\left(
\begin{array}{cc}
rr'+\mu(u\otimes v') & ru'+us' \\
vr'+sv' & \nu(v\otimes u')+ss' \\
\end{array}
\right)
$$ is a unital ring, which is equivalent to say that the maps $\mu$ and $\nu$ satisfy the
following two associativity conditions:
\begin{align}\label{associ-cond}
u\nu(v\otimes u')=\mu(u\otimes v)u'\quad \text{and}\quad v\mu(u\otimes v')=\nu(v\otimes u)v'.
\end{align}

We say that this context is \emph{strict} if the maps $\mu$ and $\nu$ are isomorphisms and, in this case,
the categories $_R\text{Mod}$ and $_S\text{Mod}$ are equivalent via the mutually inverse
equivalences $V\otimes_R-:\!_R\text{Mod}\to\!_S\text{Mod}$ and $U\otimes_S-:\!_S\text{Mod}\to\!_R\text{Mod}$. When it happens we also say that
the rings $R$ and $S$ are \emph{Morita equivalent}. If $_RU$ is faithfully projective,
that is, $_RU$ is faithful, projective and finitely generated, then it is enough the surjectivity  of  $\mu$ and $\nu$
in order to ensure the strictness of the above context. Similar statements equally hold for right module categories.

\subsection{About $B_x\star_{\bt_{(x)}}\G(x)$ and $A\star_\bt\G$}
In order to relate the partial skew groupoid ring $A\star_\bt\G$ and the partial skew group ring $B_x\star_{\bt_{(x)}}\G(x)$ we firstly recall that
\begin{align}\label{for:ideal-action}
B_g=\gamma_{\ta_{t(g)}}(I_{g_x})\quad\text{and}\quad \beta_g(\gamma_{\ta_{s(g)}}(a))=\gamma_{\ta_{t(g)}}(\gamma_{g_x}(a)), \,\,\,\text{for all}\,\,\,g\in\G,\,\, a\in I_{g\m_x}.
\end{align}

\begin{lem1}\label{lem-mc-skew} Let $R:=A\star_\bt\G$ and $S:=B_x\star_{\bt_{(x)}}\G(x)$. Then
	\begin{enumerate}[\rm (i)]
		\item  $R1_S=\sum_{g\in\mathcal{S}_x} B_g\delta_g$,\vspace*{.1cm}
		\item  $1_SR=\sum_{g\in\mathcal{T}_x} B_g\delta_g$,\vspace*{.1cm}
		\item  $1_SR1_S=S$ and $R1_SR=R$.
	\end{enumerate}
\end{lem1}
\begin{proof} Let $g\in \mathcal{S}_x$ and $a\in B_g$. By \eqref{for:ideal-action}, there exists $a'\in I_{g_x}$ such that $a=\gamma_{\ta_{t(g)}}(a')$ and whence
	\begin{align*}
	(a\delta_g)1_S&=(\gamma_{\ta_{t(g)}}(a')\delta_g)(1_x\delta_x)\\
	&=\gamma_{\ta_{t(g)}}(a')\beta_g(\gamma_{\ta_{s(g)}}(1_{g\m_x}))\delta_{g}\\
	&\overset{\mathclap{\eqref{for:ideal-action}}}{=}\gamma_{\ta_{t(g)}}(a')\gamma_{\ta_{t(g)}}(\gamma_{g_x}(1_{g\m_x}))\delta_{g}\\
	&=\gamma_{\ta_{t(g)}}(a'1_{g_x})\delta_{g}\\
	&=a\delta_g.
	\end{align*}
	For $g\notin \mathcal{S}_x$ and $a\in B_g$, we have $(a\delta_g)1_S=0$. Consequently, if $r=\sum_{g\in \G}a_g\delta_g\in R$ then $r1_S=\sum_{g\in \mathcal{S}_x}a_g\delta_g$ and (i) follows. The proof of (ii) is similar.

\vu
	
	For (iii), it is clear that $R1_SR\subset R$. For the reverse inclusion notice that
	\[\gamma_{\ta_{t(g)}}(a')\delta_g=\gamma_{\ta_{t(g)}}(a')\gamma_{\ta_{t(g)}}(1_{g_x})\delta_g=(\gamma_{\ta_{t(g)}}(a')\delta_{\ta_{t(g)}})(1_{g_x}\delta_{\ta\m_{t(g)}g}),\]
	for all $g\in \G$ and $a'\in I_{g_x}$. By (ii), $1_{g_x}\delta_{\ta\m_{t(g)}g}\in 1_SR$. Hence $\gamma_{\ta_{t(g)}}(a')\delta_g\in R1_SR$.
\end{proof}

\vu

Fix the $(R,S)$-bimodule $U:=R1_S$ and the $(S,R)$-bimodule $V:=1_SR$. Define also the bimodule map $\mu: U\otimes_S V\to R$ (resp. $\nu:V\otimes_R U\to S$) given by $u\otimes v\mapsto uv$ (resp. $v\otimes u\mapsto vu$), for all $u\in U$, $v\in V$. With this notation we have the following result.

\begin{teo1}\label{teo:mc-skew} The sixtuple $(A\star_\bt\G, B_x\star_{\bt_{(x)}}\G(x),U, V,\mu,\nu)$ is a Morita context and $\mu,\nu$ are surjective.
\end{teo1}
\begin{proof}
	It is straightforward to check that $\mu,\nu$ satisfy \eqref{associ-cond}. By Lemma \ref{lem-mc-skew} (iii), $\mu$ and $\nu$ are surjective maps.
\end{proof}

\subsection{About $A^\bt$ and $A\star_\bt\G$}

According to \cite[$\S$ 4]{BP} $A$ is a  $(A^{\beta}, A\star_{\beta}\G)$-bimodule and  a $(A\star_{\beta}\G, A^{\beta})$-bimodule. The left and the right actions of $A^\bt$ on $A$ are given by the multiplication of $A$, and the left (resp., right) action of $A\star_{\beta}\G$ on $A$ is given by $a\delta_g\cdot b=a\beta_{g}(b1_{g\m})$ (resp., $b\cdot a\delta_g=\beta_{g\m}(ba)$), for all $a,b \in A$ and $g\in \G$.

\vu

 Moreover, the maps
\[\Gamma:A\otimes_{A\star_{\beta} \G} A\to A^{\beta}, \qquad a\otimes b\mapsto t_{\beta}(ab)\]
and
\[\Gamma':A\otimes_{A^{\beta}} A\to A\star_{\beta} \G, \qquad a\otimes b\mapsto \sum_{g\in \G}a\beta_g(b1_{g\m})\delta_g,\]
are bimodule morphisms. By  \cite[Proposition 4.4]{BP}, the sixtuple $(A\star_{\beta} \G, A^{\beta}, A, A,\Gamma,\Gamma')$ is a Morita context. Furthermore, the sixtuple  $(B_x\star_{\beta_{(x)}} \G(x), B_x^{\beta_{(x)}}, B_x, B_x,\Gamma_x,\Gamma'_x)$ is also a Morita context constructed in a similar way as the previous one for the partial group action of $\G(x)$ on $B_x$ (see \cite[Theorem 1.5]{AL}).

\begin{cor1} \label{61}
$\Gamma$ is surjective if and only if $\Gamma_x$ is surjective.
\end{cor1}
\begin{proof}
Since $A$ (resp. $B_x$) is a unital ring, it is easy to see that $\Gamma$ (resp. $\Gamma_x$) is surjective if and only if $t_{\bt}$ (resp. $t_{\bt_{(x)}}$) is surjective. Hence, the result follows from Theorem \ref{t54}.	
\end{proof}

\begin{lem1} \label{62}
$\Gamma'$ is surjective if and only if $\Gamma'_{x}$ is surjective.
\end{lem1}
\begin{proof}
($\Rightarrow$) Let $b\delta_h\in B_x\star_{\bt_{(x)}}\G(x)$ with $b\in B_h\subseteq B_x$.  Consider $b\delta_h$ as in element of $A\star_\bt\G$ and the ring isomorphism
$\varphi:A\star_\bt\G\to (B\star_{\bt^\star}\G_0^2)\star_ \theta\G(x),\,\ a\delta_g\mapsto a\delta_{(s(g),t(g))}\delta_{g_x}$, given in Theorem \ref{la grande finale}. Since $\Gamma'$ is surjective, there exist $a_i,b_i\in A$, $1\leq i\leq r$, such that
$$\varphi\circ\Gamma'\left(\sum_{1\leq i\leq r}a_i\otimes b_i\right)=\varphi(b\delta_h)=b\delta_{(x,x)}\delta_h.$$
We can assume that $a_i=\somaz\gamma_{\ta_z}(a'_{i,z})$ and $b_i=\somaz\gamma_{\ta_z}(b'_{i.z})$ with $a_i,b_i\in B_x$, for all $1\leq i\leq r$. Then \vu
\begin{align*}
\varphi\circ\Gamma'\left(\sum_{1\leq i\leq r}a_i\otimes b_i\right)&=\somai\somag a_i\bt_g(b_i1_{g\m})\delta_{(s(g),t(g))}\delta_{g_x}\\
&=\somai\somag a_i\gamma_{\ta_{t(g)}}\gamma_{g_x}\gamma_{\ta\m_{s(g)}}(\gamma_{\ta_{s(g)}}(b'_{i,s(g)}1_{g\m_x}))\delta_{(s(g),t(g))}\delta_{g_x}\\
&=\somai\somag a_i\gamma_{\ta_{t(g)}}\gamma_{g_x}(b'_{i,s(g)}1_{g\m_x})\delta_{(s(g),t(g))}\delta_{g_x}\\
&=\somai\somag\gamma_{\ta_{t(g)}}(a'_{i,t(g)}\gamma_{g_x}(b'_{i,s(g)}1_{g\m_x}))\delta_{(s(g),t(g))}\delta_{g_x}.\
\end{align*} \vu
Hence, $$b\delta_{(x,x)}\delta_h=\somai\somag\gamma_{\ta_{t(g)}}(a'_{i,t(g)}\gamma_{g_x}(b'_{i,s(g)}1_{g\m_x}))\delta_{s(g),t(g))}\delta_{g_x}$$ which implies that
$$\somai a'_{i,x}\gamma_{h}(b'_{i,x}1_{h\m})=b\,\,\text{ and}\,\, \somai a'_{i,x}\gamma_l(b'_{i,x}1_{l\m})=0$$ if $l\neq h$, $l\in \G$.
Thus $\Gamma'_x(\somai a'_{i,x}\otimes b'_{i,x})=b\delta_h$ and $\Gamma'_x$ is surjective.

\vd

($\Leftarrow$) Conversely it is enough to check that for any element of the form $v= \gamma_{\ta_z}(a)\delta_{(y,z)}\delta_h$ in $(A\star_{\bt^\star}\G_0^2)\star_ \theta\G(x)$,
with $a\in A_h$ and $y,z\in\G_0$, there exists  $w\in A\otimes A$ such that $\varphi\circ\Gamma'(w)=v$. Observe that $a\delta_h\in B_x\star\bt_{(x)}\G(x)$ and whence there exist  $a_{x,i}, b_{x,i}\in B_x$ such that $a\delta_h =\Gamma'_x(\somai a_{x,i}\otimes b_{x,i})=\somal\somai a_{x,i}\gamma_l(b_{x,i}1_{l\m})\delta_l$. Thus
\begin{align}\label{for-use}
 \somai a_{x,i}\gamma_h(b_{x,i}1_{h\m})=a \ \ \text{and} \ \ \somai a_{x,i}\gamma_l(b_{x,i}1_{l\m})=0 \ \ \text{if} \ l\neq h,\ l\in \G(x).
\end{align}
Now setting $a_i=\gamma_{\ta_z}(a_{x,i})$ and $b_i=\gamma_{\ta_y}(b_{x,i})$ in $A$ one has
\begin{align*}
\varphi\circ\Gamma'\left(\somai a_i\otimes b_i\right)&=\somai\somag a_i\bt_g(b_i1_{g\m})\delta_{(s(g),t(g))}\delta_{g_x}\\
&=\somai\somag\gamma_{\ta_z}(a_{x,i})\gamma_{\ta_{t(g)}}\gamma_{g_x}\gamma\m_{\ta_{s(g)}}(\gamma_{\ta_y}(b_{x,i})\gamma_{\ta_{s(g)}}(1_{g\m_x})\delta_{(s(g),t(g))}\delta_{g_x}\\
&=\somai\sum_{g\in\G(y,z)}\gamma_{\ta_z}(a_{x,i}\gamma_{g_x}(b_{x,i}1_{g\m_x}))\delta_{(y,z)}\delta_{g_x}\\
&=\sum_{g\in\G(y,z)}\gamma_{\ta_z}(\somai a_{x,i}\gamma_{g_x}(b_{x,i}1_{g\m_x}))\delta_{(y,z)}\delta_{g_x}\\
&\overset{\eqref{for-use}}{=}\gamma_{\ta_z}(a)\delta_{(y,z)}\delta_h.\
\end{align*}
\end{proof}

\begin{teo1}\label{t63}
If $A$ is a finitely generated projective (resp., $B_x$) left $A^\bt$ (resp., $B_x^{\bt_{(x)}}$)-module   then the following statements are equivalent:
\begin{enumerate}[\rm (i)]
	\item The Morita context  $(A^{\beta}, A\star_{\beta} \G, A, A,\Gamma,\Gamma')$ is strict.
	\item The Morita context  $(B_x^{\beta_{(x)}}, B_x\star_{\beta_{(x)}} \G(x),  B_x, B_x,\Gamma_x,\Gamma'_x)$ is strict.
\end{enumerate}
\end{teo1}
\begin{proof}
It follows from Lemmas \ref{61}	and \ref{62}.
\end{proof}
\subsection{Galois theory}
The notion of partial Galois extension for partial groupoid actions was introduced in \cite{BP} as a generalization of the classical one for group actions due to S.U. Chase, D.K. Harrison and A. Rosenberg in the global case \cite{CHR} and to M. Dokuchaev, M. Ferrero and the second author in the partial case \cite{DFP}. This notion in both cases is very closed related with the strictness of the Morita context associated to it as well ensured by \cite[Theorem 5.3]{BP} and \cite[Theorem 3.1]{AL}.\vu

We say that $A$ is a \emph{$\bt$-partial Galois extension} of its subring of invariants $A^\bt$ (recalling that $\G$ was assumed to be finite) if there exist elements $a_i,b_i\in A$, $1\leq i\leq r$, such that $\somai a_i\bt_g(b_i1_{g\m})=\delta_{y,g}1_y$, for all $y\in\G_0$, where $\delta_{y,g}$ denotes the Kronecker symbol. This is equivalent to say that the map
$\Gamma'$ defined in the previous subsection is a ring isomorphism (cf. assertion (vi) in  \cite[Theorem 5.3]{BP}). Furthermore, in this case $A$ is also projective and finitely generated as a right $A^\bt$-module. (cf. assertion (ii) in  \cite[Theorem 5.3]{BP}).

\vu

The following theorem shows how really closed are the Galois theory for partial connected groupoid actions and the one for partial group actions.

\begin{teo1}\label{76}
The following statements are equivalent:
\begin{enumerate}[\rm (i)]
	\item $A$ is a partial Galois extension of $A^\bt$ and $t_\bt$ is onto. \vu
	\item $(A^{\beta}, A\star_{\beta} \G, A, A,\Gamma,\Gamma')$ is strict.\vu
	\item $(B_x^{\beta_{(x)}}, B_x\star_{\beta_{(x)}} \G(x),  B_x, B_x,\Gamma_x,\Gamma'_x)$ is strict.\vu
	\item $B_x$ is a partial Galois extension of $B_x^{\bt_{(x)}}$ and $t_{\bt_{(x)}}$ is onto.
\end{enumerate}
\end{teo1}
\begin{proof}
It is an immediate consequence of  Theorem \ref{t54}, Theorem \ref{t63},  of \cite[Theorem 5.3]{BP} and   \cite[Theorem 3.1]{AL}.
\end{proof}
\section{Separability, Semisimplicity and Frobenius properties}

In this section we analyze the properties of separability, Frobenius and semisimplicity concerning to the extensions $A\subset A\star_{\beta}\G$ and $B_x\subset B_x\star_{\beta_{(x)}} \G(x)$.
\subsection{Separability and semisimplicity} A unital ring extension $R\subset S$ is called {\it separable} if the multiplication map $m:S\otimes_R S\to S$ is a splitting epimorphism of $S$-bimodules. This is equivalent to say that there exists an {\it idempotent of separability} of $S$ over $R$, i.~e. an element $x\in S\otimes_RS$ such that $sx=xs$, for all $s\in S$, and $m(x)=1_S$.
A ring extension $R\subset S$ is left (right) {\it semisimple} if any left (right) $S$-submodule $N$ of a left (right) $S$-module $M$ having an $R$-complement
in $M$, has an $S$-complement in $M$. Further details on separable and semisimple extensions can be seen in \cite{DI}, \cite{HS} and \cite{CIM}.

\begin{remark}\label{ob:sepa-semi}
	\rm{It is well-known that any separable extension is a left (right) semisimple extension (see for instance \cite{CIM}).}
\end{remark}

%
%Given a partial action $\af=(A_g,\af_g)_{g\in G}$ of a finite group $G$ on a ring $A$, we recall that the trace map $t_{\alpha}:A\to A$ is defined by $t_{\alpha}(a)=\sum_{g\in G}\af_g(a1_{g\m})$, for all $a\in A$.

 For the partial action $\bt$ of $\G$ on $A$ we associate the maps, as introduced in \cite{BPi}, $t_{y,z}:A\to A$ and $t_{z}:A\to A$ given by
\[ t_{y,z}(a)=\sum_{g\in \in \G(y,z)} \beta_g(a1_{g\m}), \quad t_z(a)=\sum_{y\in \G_0} t_{y,z}(a),\quad y,z\in \G_0,\,\,a\in A.\]

\vu

Notice that $\varphi\colon \A\to \A\star_{\bt}\G$, $\ a\mapsto \sum_{y\in \G_0}(a1_y)\delta_y$, is a monomorphism of rings and whence $\A\star_{\bt}\G$ is an extension of $A$. Moreover, $\varphi$ induces the following $(\A,\A)$-bimodule structure on $\A\star_{\bt}\G$:
\begin{align}\label{aabim}
&a\cdot (a_g\delta_g)=aa_g\delta_g,\quad\quad
(a_g\delta_g)\cdot a=a_g\bt_g(a1_{g^{-1}})\delta_g,& &g\in\G,\ a\in \A.&
\end{align}

Now we can present the main result of  this subsection.

\vu

\begin{teo1} \label{82} The following statements are equivalent:
	\begin{enumerate}[\rm (i)]
		\item there is an element $a$ in the center $C(A)$ of $A$ such that $t_{y}(a)=1_y$ for all $y\in \G_0$,
\vu
		\item there is an element $b$ in the center $C(B_x)$  of $B_x$ such that $t_{\beta_{(x)}}(b)=1_x$,
\vu
		\item $B_x\subset B_x\star_{\beta_{(x)}} \G(x)$ is separable,
\vu
		\item $A\subset A\star_{\beta}\G$ is separable.
	\end{enumerate}	
\end{teo1}

\begin{proof}
	
	\noindent (i) $\Rightarrow$ (ii) Let $a\in C(A)$ such that $t_{y}(a)=1_y$ for all $y\in \G_0$. Observe that $a=\sum_{z\in \G_0} \gamma_{\tau_z}(a'_z)$ with $a'_z\in C(B_x)$ for each $z\in \G_0$. Also,
	\begin{align*}
	t_{y,x}(a)&=\sum_{g\in \G(y,x)}\gamma_{\ta_x}\gamma_{g_x}\gamma\m_{\ta_y}(a\gamma_{\tau_y}(1_{g\m_x}))\\
	&=\sum_{g\in \G(y,x)}\gamma_{g_x}\gamma\m_{\ta_y}(\gamma_{\tau_y}(a'_y1_{g\m_x}))\\
	&=\sum_{g\in \G(y,x)}\gamma_{g_x}(a'_y1_{g\m_x})\\
	&=\sum_{h\in \G(x)}\gamma_{h}(a'_y1_{h\m})\qquad (h=g\tau_y)\\
	&=t_{\beta_{(x)}}(a'_y).
	\end{align*}
	Hence, $1_x=t_{x}(a)=\sum_{y\in \G_0}t_{\beta_{(x)}}(a'_y)=t_{\beta_{(x)}}(b)$, where $b=\sum_{y\in \G_0}a'_y\in C(B_x)$.\vu
	
	\noindent (ii) $\Rightarrow$ (i) Conversely, assume that there exists $b\in C(B_x)$ such that
	$t_{\beta_{(x)}}(b)=1_x$. Notice that $b\in C(A)$, $t_{y,z}(b)=0$ if $y\neq x$ and
	\[t_{x,z}(b)=\sum_{g\in \G(x,z)}\gamma_{\ta_z}\gamma_{g_x}(b1_{g\m_x})=\gamma_{\ta_z}(t_{\beta_{(x)}}(b)=\gamma_{\ta_z}(1_x)=1_z.\]
	Consequently, $t_z(b)=1_z$ for all $z\in \G_0$. \vu
	
	Moreover, by    \cite[Theorem 3.1]{BLP} we have that  (ii) $\Leftrightarrow$ (iii) and  \cite[Theorem 4.1]{BPi} implies  that (i) $\Leftrightarrow$ (iv).
\end{proof}

\begin{cor1} \label{83}
	If the equivalent statements of Theorem \ref{82}  hold then the ring extensions $B_x\subset B_x\star_{\beta_{(x)}} \G(x)$ and $A\subset A\star_{\beta}\G$ are semisimple.	
\end{cor1}
\begin{proof}
	It follows from Remark \ref{ob:sepa-semi} and Theorem \ref{82}.
\end{proof}

\subsection{Frobenius property}
A ring extension $R\subset S$ is called {\it Frobenius} if there exist an element $u=\sum_{i=1}^{n}s_{i1}\otimes s_{i2}\in S\otimes_R S$ and an $R$-bimodule map $\varepsilon:S\to R$ such that $x$ is $S$-central and $\sum_{i=1}^{n}\varepsilon(s_{i1})s_{i2}=\sum_{i=1}^{n}s_{i1}\varepsilon(s_{i2})=1$.  More details on Frobenius extension can be seen, for example, in \cite{CIM}.

\begin{teo1}\label{84} The extensions $A\subset A\star_{\beta}\G$ and $B_x\subset B_x\star_{\beta_{(x)}} \G(x)$ are Frobenius.
\end{teo1}
\begin{proof} Consider $u=\sum_{g\in \G}1_g\delta_g\otimes 1_{g^{-1}}\delta_{g^{-1}}\in A\star_{\beta}\G\otimes_{A} A\star_{\beta}\G $ and note that
	\begin{align*}
	(a_l\delta_l) u&=\sum_{g\in \G} a_l\delta_l1_g\delta_g\otimes 1_{g^{-1}}\delta_{g^{-1}} \\
	&=\sum_{g\in \mathcal{T}_{s(l)}} a_l\beta_l(1_g1_{l^{-1}})\delta_{lg}\otimes 1_{g^{-1}}\delta_{g^{-1}}\\
	&=\sum_{h\in \mathcal{T}_{t(l)}} a_l1_h\delta_{h}\otimes 1_{h^{-1}l}\delta_{h^{-1}l} \qquad\,\, (h=lg)\\
	&\overset{\mathclap{\eqref{aabim}}}{=}\sum_{h\in \mathcal{T}_{t(l)}} 1_h\delta_{h}\cdot \beta_{h\m}(a_l1_h)\otimes 1_{h^{-1}l}\delta_{h^{-1}l}  \\
	&=\sum_{h\in \mathcal{T}_{t(l)}} 1_h\delta_{h} \otimes \beta_{h\m}(a_l1_h)\delta_{h^{-1}l}\\
	&=\sum_{h\in \mathcal{T}_{t(l)}} 1_h\delta_{h} \otimes 1_{h\m}\delta_{h\m} a_l\delta_{l}\\
	&=u(a_l\delta_l),
	\end{align*}	
	for all $l\in \G$ and $a_l\in B_l$. Hence, $u$ is $A\star_{\beta}\G$-central. Note that $\varepsilon:A\star_{\beta}\G\to A$	given by \[\varepsilon(a_g\delta_g)=\delta_{g,\G_0}a_g=\left\{\begin{array}{lc}
	a_g, 	& \text{ if } g\in \G_0\\
	0,& \text{ otherwise}
	\end{array}
	\right. \] is an $R$-bimodule map. Moreover, \[\sum_{g\in\G}\varepsilon(1_g\delta_g)1_{g\m}\delta_{g\m}=\sum_{g\in\G}1_g\delta_g1_{g\m}\varepsilon(\delta_{g\m})=\sum_{y\in \G_0}1_y\delta_y=1.\]
	Thus $A\subset A\star_{\beta}\G$ is a Frobenius extension. The second assertion is an immediate consequence of   \cite[Theorem 3.6]{BLP}.
\end{proof}	

\end{document}